\documentclass[11pt, a4paper, oneside]{article}
\usepackage[body={15.5cm, 21.0cm},left=2cm,right=2cm]{geometry}
\usepackage{natbib}
\usepackage{amsmath, amssymb}
\usepackage{charter}
\usepackage{graphicx}
\usepackage{setspace}	
\usepackage{pifont}
\usepackage{color}
\usepackage{hyperref}
\usepackage{fancyhdr}
\newtheorem{thm}{Theorem}
\newtheorem{lem}[thm]{Lemma}
\newtheorem{rem}[thm]{Remark}

\newtheorem{prop}[thm]{Proposition}
\newtheorem{ex}[thm]{Example}

\newenvironment{pf}{\noindent{\bf Proof:}}{\hfill \ding{111} \\}
\newenvironment{pfofThm}{\noindent{\bf Proof of Theorem }}{\hfill \ding{111} \\}

\newcommand{\ie}{i.e., }

\newcommand{\id} {\ensuremath{\displaystyle{\mathop {=} ^d}}}

\newcommand{\field}[1]{\mathbb{#1}}
\newcommand{\real}{\ensuremath{{\field{R}}}}
\newcommand{\mc}[1]{{\ensuremath{\mathcal{#1}}}}

\newcommand{\sumab}[2]{\ensuremath{\sum\limits_{#1}^{#2}}}
\newcommand{\intab}[2]{\ensuremath{\int_{#1}^{#2}}}
\newcommand{\intinf}[1]{\ensuremath{\int_{#1}^{\infty}}}
\newcommand{\intunit}{\ensuremath{\int_{0}^{1}}}

\newcommand{\arrowf}[1]{\ensuremath{\displaystyle {\mathop {\longrightarrow}_{#1 \rightarrow \infty}\,}}}
\newcommand{\arrowd}[1]{\ensuremath{\displaystyle {\mathop {\longrightarrow}_{#1}\,}}}
\newcommand{\limit}[1]{\ensuremath{\displaystyle {\lim_{#1 \rightarrow{\infty}}}}}
\newcommand{\suprem}[1]{\ensuremath{\displaystyle {\sup_{#1}}}}

\newcommand{\conv}[1]{\ensuremath{\, \displaystyle {\mathop {\longrightarrow}_{n \rightarrow \infty} ^{#1}}}\, }

\newcommand{\ndivk}{\frac{n}{k}}

\title{Estimation of the finite right endpoint in the Gumbel domain}

\author{{Isabel Fraga Alves}\\  {\small CEAUL and DEIO}\\ {\small FCUL, University of Lisbon} \and {Cl\'{a}udia Neves}\\ {\small CEAUL and University of Aveiro}}

\date{}

\begin{document}

\maketitle

\strut

\abstract{A simple estimator for the finite right endpoint of a distribution function in the Gumbel max-domain of attraction is proposed. Large sample properties such as consistency and the asymptotic distribution are derived. A simulation study is also presented.}




\section{Introduction}
\label{SecIntro}

Let $X_{n,n}\geq X_{n-1,n}\geq \ldots \geq X_{1,n}$ be the order statistics from the sample $X_1, X_2,\ldots,X_n$ of i.i.d. random variables with common (unknown) distribution function $F$. Let $x^F$ denote the right endpoint of $F$. We shall assume that the distribution function $F$ has a finite right endpoint, i.e. $x^F:= \sup\{x:\, F(x)<1\} \in \real$.

The fundamental result for extreme value theory is due in vary degrees of generality to \citet{FisherTippett:28}, \citet{G:43}, \citet{deHaan:70} and \citet{BdH:74}. The extreme value theorem (or extremal types theorem) surprisingly restricts the class of all possible limiting distribution functions to only three different types, while the induced domains of attraction embrace a great variety of distribution functions. This is particularly true in the case of the Gumbel domain of attraction. In other words, if there exist constants $a_n>0$, $b_n\in \real $ such that
\begin{equation}\label{MaxDom}
	\limit{n} F^n(a_n\,x+b_n)= G(x),
\end{equation}
for all $x$, $G$ non-degenerate, then $G$ must be only one of the following:
\begin{eqnarray*}
	\Psi_\alpha(x)&=&\exp\{-(-x)^\alpha\},\;x<0,\quad \alpha>0,\\
	\Lambda(x) &=& \exp\{-\exp(-x)\},\quad x\in\real,\\
	\Phi_\alpha(x)&=&\exp\{-x^{-\alpha}\},\;x>0,\quad \alpha>0.
\end{eqnarray*}
Redefining the constants $a_n>0$ and $b_n\in \real$, these can in turn be nested in a one-parameter family of distributions, the Generalized Extreme Value (GEV) distribution with distribution function
\begin{equation*}
	G_{\gamma}(x):= \exp\{-(1+\gamma x)^{-1/\gamma}\}, \; 1+\gamma x>0,\;
\gamma \in \mathbb{R}.
\end{equation*}
We then say that $F$ is in the (max-)domain of attraction of $G_\gamma$ and use the notation $F\in \mc{D}_{M}(G_{\gamma})$. For $\gamma<0$, $\gamma=0$ and $\gamma>0$, the GEV distribution function reduces again to Weibull, Gumbel and Fr\'echet distribution functions, respectively. An equivalent extreme value condition allows the limit relation in \eqref{MaxDom} to
run over the real linear \citep[cf. Theorem 1.1.6][]{deHaanFerreira:06}: $F\in \mathcal{D}_{M}(G_{\gamma})$ if and only if
\begin{equation}\label{DOA}
	\limit{t} t\bigl( 1-F(a(t)\,x+b(t)\bigr)=(1+\gamma\, x)^{-1/\gamma},
\end{equation}
for all $x$ such that $1+\gamma x>0$, $a(t):= a_{[t]}$ and $b(t):= b_{[t]}$, with $[t]$ denoting the integer part of $t$. The extreme value index $\gamma$ determines vary degrees of tail heaviness. If $F\in \mathcal{D}_{M}(G_{\gamma})$ with $\gamma>0$, then the distribution function $F$ is heavy-tailed, i.e., $F$ has a power-law decaying tail with infinite right endpoint. On the opposite end, $\gamma <0$ refers to short tails which must have finite right endpoint. The Gumbel domain of attraction $\mc{D}_{M}(G_0)$ renders a great variety of distributions, ranging from light-tailed distributions such as the Normal distribution, the exponential distribution, to moderately heavy distributions such as the Lognormal. All the just mentioned distributions have an infinite right endpoint but a finite endpoint is also possible in the Gumbel domain. We shall give several examples. Distribution functions of this sort, i.e. light-tailed distributions with finite endpoint, but not so light that they are still included in the Gumbel domain, have been in great demand as feasible distributions underlying real life phenomena. A striking example is the extreme value analysis by \citet{EM:08} of the best marks in Athletics, aiming at assessing the ultimate records for several events.
For instance, Table 3 in \citet{EM:08} has several missing values for the estimates of the endpoint which are due to an estimated extreme value index $\gamma$ near zero. An attempt to fulfill these blank spaces  with an appropriate framework for inference in the Gumbel domain has been provided by \citet{FAdHN:10}, althought from the strict view point of application to the Long Jump data set used in \citet{EM:08}. The tentative estimator proposed by \citet{FAdHN:10} is virtually the same as the one introduced in the present paper. The novelty here is in the development of a simple closed-from expression for the previous statistic. Hence, the problem of estimating the right endpoint $x^F$ of a distribution function lying in the Gumbel extremal domain of attraction is now tackled by the semi-parametric statistic
\begin{equation*}
	X_{n,n}+X_{n-k,n}- \frac{1}{\log 2}\sumab{i=0}{k-1}\log \Bigl(1+\frac{1}{k+i}\Bigr) X_{n-k-i,n},
\end{equation*}
or in a more compact form, by
\begin{equation}\label{Combi}
	\hat{x}^F:= X_{n,n}+  \sumab{i=0}{k-1}a_{i,k} \Bigl( X_{n-k,n}-X_{n-k-i,n}\Bigr),
\end{equation}
where $a_{i,k}:=\bigl( \log(k+i+1)-\log (k+i)\bigr)/\log 2$, such that $\sum_{i=0}^{k-1} a_{i,k}=1 $. Here and throughout this paper, the number $k$ is assumed intermediate, that is, $k$ is in fact a sequence of positive integers going to infinity as $n\rightarrow \infty$ but at a much slower rate than $n$. More formally, we are assuming that $\hat{x}^F$ is a functional of the top observations of the original sample, which relies on an intermediate sequence $k = k_n$ such that
\begin{equation*}
k_{n}\rightarrow \infty, \quad k_n=o(n), \quad \mbox{ as } n\rightarrow \infty.
\end{equation*}
From the non-negativeness of the weighted spacings in the sum \eqref{Combi}, we clearly see that the now proposed estimator is greater than $X_{n,n}$ with probability one. This constitutes a crucial advantage in comparison with the usual semi-parametric estimators for the right endpoint of a distribution function in the Weibull domain of attraction (i.e. with $\gamma<0$). We refer to \citet{Hall:82}, \citet{Falk:95}, \citet{HallWang:99} and to \citet{deHaanFerreira:06} and references therein. To the best of our knowledge, none of these estimators have ensured so far the extrapolation beyond the sample range, meaning that we can encounter in practice estimates for the endpoint that are smaller than the observed sample maximum. There have been, however, some developments of the most well-known endpoint estimators connected with $\gamma<0$ in the sense of bias reduction and/or correction. \cite{LiPeng:09} and \citet{Caietal:12} are two of the most recent works in this respect. In fact, the problem of estimating $x^F$ still gathers a great interest nowadays. Recently, \citet{Guillouetal:12} devised an endpoint estimator from the high-order moments pertaining to a distribution attached with $\gamma<0$; \citet{LiPeng:12} proposed a bootstrap estimator for the endpoint evolving from the one by \citet{Hall:82} in case $\gamma \in (-1/2,0)$. The present paper deliberately addresses the class of distribution functions belonging to the Gumbel domain of attraction, for which no specific inference has yet been provided in the context of estimation of the right endpoint $x^F < \infty$. The appropriate framework for the latter shall be developed in Section \ref{SecFrame}.

The remainder of the paper is as follows. The rationale behind the proposal of the new estimator for the right endpoint is expounded in Section \ref{SecStats}. Large sample properties of this estimator, namely consistency and asymptotic distribution, are worked out in Section \ref{SecAsympt} by taking advantage of this form of separability between the maximum and the sum of higher order statistics. In order to perform asymptotics, we require some basic conditions in the context of the theory of regular variation. These are laid out in the next section (Section \ref{SecFrame}). Finally, in Section \ref{SecSims} we gather some simulation results taken as key examples.


\section{Framework}
\label{SecFrame}
Let $F$ be a distribution function (d.f.) with right endpoint $x^F$,
\begin{equation*}
x^F:= \sup\{x:\, F(x)<1\}.
\end{equation*}
Suppose $F$ belongs to the domain of attraction of the Generalized Extreme Value distribution (GEV) with d.f. $G_{\gamma}$, that is $F$ satisfies the following extreme value condition
\begin{equation}\label{EVcond}
\displaystyle {\lim_{x \uparrow x^F}}\, \frac{1-F(t+x\,f(t))}{1-F(t)}= (1+\gamma\, x)^{-1/\gamma},
\end{equation}
for all $x\in \real$ such that $1+\gamma\,x>0$, with a suitable positive function $f$ \citep[equivalent condition to \eqref{DOA}, see Theorem 1.1.6 of][]{deHaanFerreira:06}.

For the most interesting case of $\gamma=0$ the limit in \eqref{EVcond} reads as $e^{-x}$. In this case $f>0$ can be defined as follows
\begin{equation}\label{fFunct}
f(t):= \intab{t}{x^F}\frac{1-F(x)}{1-F(t)}\, dx= E[X-t|X>t]
\end{equation}
 \citep[cf. Theorem 1.2.5 of][]{deHaanFerreira:06}, then $f$ is the so called Mean Excess Function.

Now let $U$ be the (generalized) inverse function of $1/(1-F)$.
If $F$ satisfies \eqref{EVcond} with $\gamma=0$ then we can assume there exists a positive function $a_0$ such that, for all $x>0$,
\begin{equation}\label{PiofU}
\limit{t}\frac{U(tx)-U(t)}{a_0(t)}= \log x.
\end{equation}
Hence $U$ belongs to the class $\Pi$  \citep[see Definition B.2.4 of][]{deHaanFerreira:06} and $a_0$ is a measurable function such that $\limit{t}  a_0(tx)/a_0(t)=1$ for
$x>0$. Then we say that $a_0$ is a slowly varying function and use the notation $a_0 \in RV_0$ \citep[see Theorem B.2.7 of][]{deHaanFerreira:06}. Moreover,  the functions $a_0$ and $f$ (introduced in \eqref{PiofU} and \eqref{EVcond}, respectively) are related to each other by  $a_0=f \circ U$\citep[see Theorem B.2.21 of][]{deHaanFerreira:06}. Throughout we shall use the notation $U\in \Pi (a_0)$ in order to put some emphasis on the auxiliary function $a_0$. We have the following result:

\begin{lem} \label{LemStem}\mbox{ }
\begin{enumerate}
  \item  Suppose $U\in \Pi(a)$. For any $\varepsilon>0$ there exists $t_0=t_0(\varepsilon)$ such that, for $st\geq t_0$,
\begin{equation*}
\biggl| \frac{a(st)}{a(t)}-1\biggr| \leq \varepsilon\, \max(s^{\varepsilon},s^{-\varepsilon}).
\end{equation*}
  \item Suppose $a>0$ is a slowly varying function,  integrable over finite intervals of $\real^+$ such that
  \begin{equation*}
	\intinf{t} a(s) \, \frac{ds}{s} <\infty.
  \end{equation*}
for every $t>0$. Then $a(t) \rightarrow 0\,$, as $t\rightarrow \infty$, and
\begin{equation*}
	\limit{t}\intinf{t}\frac{a(s)}{a(t)}\,\frac{ds}{s}= \infty.
\end{equation*}
\end{enumerate}
\end{lem}
\begin{pf}
Part 1. of the Lemma comes from \citep{Drees:98} \cite[cf. Proposition B.1.10 of][]{deHaanFerreira:06}. The second part follows from Karamata's theorem for regularly varying functions \cite[cf. Theorem B.1.5 of][]{deHaanFerreira:06}.
\end{pf}

The relationship between conditions imposed on the auxiliary function $a$ (i.e. two conditions in 2. of Lemma \ref{LemStem}) and the tail quantile function $U$, for which $x^F:= U(\infty)= \lim_{t \rightarrow \infty} U(t)$ exists finite, is given by
\begin{equation}\label{CharU}
U(t) = c+\intab{1}{t} a(s)\, \frac{ds}{s} + o\bigl(a(t)\bigr),
\end{equation}
$c\in \real$ \citep[cf. Theorem B.2.12 and Proposition B.2.15 (3.) of][]{deHaanFerreira:06}. In this development, the following holds:
\begin{equation}\label{MainRel}
U(\infty)-U(t)= \intinf{t}a(s)\, \frac{ds}{s} + o\bigl(a(t)\bigr), \quad t\rightarrow \infty,
\end{equation}
which is our main assumption eventually. Moreover, \eqref{MainRel} implies that $U\in \Pi(a)$ and $a(t)\sim a_0(t)$, as $t\rightarrow \infty$, with $a_0$ the auxiliary function in \eqref{PiofU}.

We can obtain from \eqref{PiofU} with $a_0$ replaced by $a$ (i.e. $U\in\Pi(a)$) yet another limiting relation now involving integration of $U$ and $a$: applying Cauchy's rule once, we obtain
\begin{equation}\label{CauchyAppl}
\limit{t} \, \frac{\intab{0}{1/t}\Bigl(U\bigl(\frac{x}{s} \bigr)-U\bigl(\frac{1}{s} \bigr)\Bigr)\,\frac{ds}{s}}{\intab{0}{1/t}a\bigl(\frac{1}{s}\bigr)\,\frac{ds}{s}}= \limit{t}\, \frac{\bigl(U(tx)-U(t)\bigr))/t}{a(t)/t},
\end{equation}
then for arbitrary positive $x$, the $\Pi$-variation of $U$ ascertains that $\log x$ is the limit above, i.e.
\begin{equation}\label{PiInt}
\limit{t}\frac{\intinf{tx}U(s)\,\frac{ds}{s}-\intinf{t}U(s)\,\frac{ds}{s}}{\intinf{t}a(s)\,\frac{ds}{s}}= \log x,
\end{equation}
for all $x>0$. Hence $\intinf{t}U(s)\,ds/s$ is also $\Pi$-varying with auxiliary function
\begin{equation}\label{qDef}
q(t):=\intinf{t}a(s)\,\frac{ds}{s}.
\end{equation}
In the usual notation, $\intinf{t}U(s)\,ds/s \in \Pi(q)$. Then $q$ is slowly varying while relation \eqref{MainRel} entails that $q(t)\rightarrow 0$ as $t \rightarrow \infty$.

\section{Statistics}
\label{SecStats}

Let $X_1, X_2, \ldots, X_n$ be a random sample of size $n$ from the underlying distribution function $F$ with finite right endpoint $x^F$. Let  $X_{1,n}\leq X_{2,n}\leq \ldots\leq X_{n,n}$ be the corresponding order statistics. We
introduce the estimator $\hat{q}(n/k)$ for the auxiliary function defined above, i.e.
\begin{equation}\label{qAlt}
q(t)= \intinf{1} a(st)\,\frac{ds}{s}= \intunit a\bigl(\frac{t}{s}\bigr)\,\frac{ds}{s}
\end{equation}
evaluated at $t=n/k$. This estimator has the property that, as $n\rightarrow \infty$,
$k=k(n) \rightarrow \infty$ and $k(n)/n \rightarrow 0$ (provided some suitable yet
mild restrictions involving the second order refinement of $\intinf{t}U(s)/s\,ds$),
\begin{equation*}
    \frac{q\bigl(\ndivk \bigr)}{a\bigl(\ndivk \bigr)}\biggl( \frac{\hat{q}\bigl(\ndivk \bigr)}{q\bigl(\ndivk \bigr)}-1\biggr) \conv{d}
    N,
\end{equation*}
where $N$ is a non-generate random variable. Several estimators for the right endpoint $x^F= U(\infty)<\infty$ can be readily devised from \eqref{MainRel}, in the sense that  these might evolve from
\begin{equation}\label{EstEndPoint}
\hat{x}^F= \widehat{U}\bigl(\ndivk \bigr)+ \hat{q}\bigl(\frac{n}{k}\bigr)= X_{n-k,n}+ \hat{q}\bigl(\frac{n}{k}\bigr),
\end{equation}
which also enables the estimates yields to carry analogous large sample properties to $\hat{q}(n/k)$.
In particular, relation \eqref{PiInt} at $x=1/2$ together with \eqref{qAlt} at $t=n/k$ prompts the following approximation for large enough $n$:
\begin{equation*}
\intunit \Bigl(U\bigl( \frac{n}{2ks}\bigr)-U\bigl(\frac{n}{ks}\bigr)  \Bigr)\, \frac{ds}{s} \approx q\bigl( \ndivk\bigr)(-\log 2 ).
\end{equation*}
Our proposal for estimating $q(n/k)$ thus arises quite naturally from the corresponding empirical counterparts (i.e. $\hat{U}\bigl(n/(\theta ks)\bigr)=X_{n-[\theta ks],n}$, $s \in (0,1]$, $\theta=1,2$):
\begin{equation}\label{Est1}
\hat{q}\bigl( \ndivk\bigr):= -\frac{1}{\log 2} \intunit \bigl( X_{n-[2ks],n}-X_{n-[ks],n}\bigr)\, \frac{ds}{s} .
\end{equation}
A certain amount of simple calculations yields the following alternative expression for $\hat{q}$:
\begin{equation}\label{Est1Alt}
\hat{q}\bigl( \ndivk\bigr)= X_{n,n}+ \frac{1}{\log 2}\sumab{i=0}{k-1}\log \Bigl(\frac{k+i}{k+i+1}\Bigr) X_{n-k-i,n}.
\end{equation}
Combining  \eqref{EstEndPoint} with  \eqref{Est1Alt} we are led to the estimator  for the right endpoint
\begin{equation}\label{EstxF}
\hat{x}^F:= X_{n-k,n} +X_{n,n}+ \frac{1}{\log 2}\sumab{i=0}{k-1}\log \Bigl(\frac{k+i}{k+i+1}\Bigr) X_{n-k-i,n}.
\end{equation}

We note that, after rearranging components, it is possible to express $\hat{x}^F$ as the maximum $X_{n,n}$ added by some  weighted mean of non-negative summands as follows:
\begin{eqnarray*}
	\hat{x}^F&=&  X_{n,n}+  \sumab{i=0}{k-1}a_{i,k} \Bigl( X_{n-k,n}-X_{n-k-i,n}\Bigr),\\
	a_{i,k}&:=&-\frac{1}{\log 2}\bigl( \log (k+i)-\log(k+i+1)\bigr)>0,
\end{eqnarray*}
all $i=0,1,\ldots$ and $k\in \field{N}$. We can easily see that $a_{i,k}$ are such that
\begin{equation*}
	\sumab{i=0}{k-1}a_{i,k}= 1.
\end{equation*}

\begin{rem}We emphasize that the now proposed estimator for the right endpoint returns values always larger than $x_{n,n}$. This constitutes a major advantage in comparison to  the available semi-parametric estimators for the endpoint  in the case of Weibull domain of attraction, for which the extrapolation beyond the sample range is not guaranteed. This inadequacy of the existing estimators often leads to some disappointing results in practical applications, with estimates-yields that may be lower than the observed maximum in the data.
\end{rem}
\section{Asymptotic results}
\label{SecAsympt}

Throughout this section we shall bear in mind that $\hat{x}^F$ rests clearly on two building blocks: the high random threshold $X_{n-k,n}$ and $\hat{q}(n/k)$ defined in \eqref{Est1Alt}.
We shall handle $\hat{q}(n/k)$ first. The proof for consistency of the estimator $\hat{q}(n/k)$ defined in \eqref{Est1} is supported on the assertion in Lemma 2.4.10 of \citet{deHaanFerreira:06}. The asymptotic distribution of  $\hat{q}(n/k)$ is attained under a second order limit  regarding the main conditions of (extended) regular variation provided in Section \ref{SecFrame}, by taking advantage of its inherent separability between the maximum and the sum of other high-order statistics. Then the two main results concerning $\hat{x}^F$, comprising Theorem \ref{TheoConsxF} and Theorem \ref{TheoANxF}, arise almost directly from the previous.

Let $U_1, U_2, \ldots, U_n$ be independent and identically distributed uniform random variables on the unit interval and let $U_{1,n}\leq U_{2:n}\leq \ldots \leq U_{n:n}$ be their order statistics. Note that $U(1/U_{i})\id X_i$, $i=1,2,\ldots$.
Since $k=k(n)$ is an intermediate sequence such that $k(n)\rightarrow \infty$, $k(n)=o(n)$, as $n\rightarrow \infty$, then we can define a sequence of Brownian motions $\bigl\{W_n(s)\bigr\}_{s\geq 0}$ such that, for each $\varepsilon>0$,
\begin{equation}\label{UnifTailProc}
\suprem{\frac{1}{\theta k}\leq s\leq 1}s^{\frac{3}{2}+\varepsilon}\biggl|\sqrt{\theta k}\biggl(\frac{\theta k}{ n U_{[\theta k s]+1,n}}-\frac{1}{s}\biggr)-\frac{1}{s^2}  W_n(s) \biggr|= o_p(1),
\end{equation}
for all  $\theta \geq 1$ (cf. Lemma 2.4.10 of de Haan and Ferreira, 2006, with $\gamma=1$).\\

Let $X_1,\, X_2,\ldots $ be i.i.d random variables with the same distribution function $F$ belonging to the Gumbel domain of attraction, i.e., $F\in \mathcal{D}(G_0)$, with finite right endpoint $x^F$. In view of characterization \eqref{CharU} for $U\in \Pi(a)$, the following relation holds
\begin{equation*}
\frac{U(tx)-U(t)}{a(t)}\approx \intab{1/x}{1}\frac{a\bigl(\frac{t}{s}\bigr)}{a(t)}\, \frac{ds}{s}, \quad (t\rightarrow \infty)
\end{equation*}
for all $x>0$. Hence we obtain for sufficiently large $n$ that
\begin{equation*}
\frac{X_{n-[\theta k s],n}-U\bigl(\frac{n}{\theta k} \bigr)}{a\bigl(\frac{n}{\theta k}\bigr)}\,\id \, \frac{U\bigl(\frac{n}{\theta k}\,\frac{\theta k}{n U_{[\theta ks]+1,n}}\bigr)- U\bigl(\frac{n}{\theta k} \bigr)}{a\bigl(\frac{n}{\theta k}\bigr)}\\
	\,\approx\, \intab{\frac{n U_{[\theta ks]+1,n}}{\theta k}}{1}\frac{a\bigl(\frac{n}{\theta k}\frac{1}{x} \bigr)}{a\bigl(\frac{n}{\theta k }\bigr)}\,\frac{dx}{x}.
\end{equation*}
Now the uniform inequalities in Lemma \ref{LemStem}(1) tell us that, for any $\varepsilon>0$,
\begin{equation*}
\frac{a\bigl(\frac{n}{\theta k}\frac{1}{s}\bigr)}{a\bigl(\frac{n}{\theta k }\bigr)} = 1 \pm \varepsilon s^{-\varepsilon}, 0< s \leq 1.
\end{equation*}
Since {$U_{[\theta ks]+1,n}\in [0,1]$ and for every $s\in (0,1]$,
\begin{equation*}
\frac{n\,U_{[\theta ks]+1,n}}{\theta k}\leq \frac{n\,U_{[\theta k]+1,n}}{\theta k} \conv{P} 1,
\end{equation*}
 we get
\begin{eqnarray*}
\frac{X_{n-[\theta k s],n}-U\bigl(\frac{n}{\theta k} \bigr)}{a\bigl(\frac{n}{\theta k}\bigr)}&\id& -\log s - \log \biggl(1+\Bigl( \frac{n U_{[\theta ks]+1,n}}{\theta ks}-1\Bigr) \biggr)\pm \biggl(\Bigl( \frac{n U_{[\theta ks]+1,n}}{\theta k}\Bigr)^{-\varepsilon}-1\biggr)\\
& = & -\log s -\frac{1}{s}\Bigl( \frac{n U_{[\theta ks]+1,n}}{\theta k}-s\Bigr)\bigl( 1+o_p(1)\bigr)\pm (s^{-\varepsilon}-1)\bigl( 1+o_p(1)\bigr),
\end{eqnarray*}
with the $o_p(1)$-term tending to zero uniformly for $s \in [(\theta k)^{-1},1]$. Now we can apply Cram\'er's $\delta$-method  to relation \eqref{UnifTailProc} in order to obtain:
\begin{equation}\label{RepIntOrder}
\frac{X_{n-[\theta k s],n}-U\bigl(\frac{n}{\theta k} \bigr)}{a\bigl(\frac{n}{\theta k}\bigr)}= -\log s +\frac{1}{\sqrt{\theta k}}\bigl(s^{-1}W_n(s) + o_p(s^{-1/2-\varepsilon})\bigr)\pm(s^{-\varepsilon}-1)\bigl( 1+o_p(1)\bigr),
\end{equation}
as $n\rightarrow \infty$, uniformly for $(\theta k)^{-1}  \leq s \leq 1$, $\theta \geq 1$. We now consider the normalized difference between a sample intermediate quantile and corresponding theoretical quantile and denote it by $R_{\theta}(s)$, i.e.
\begin{eqnarray}
\label{RTheta} R_{\theta}(s)&:= &\frac{X_{n-[\theta k s],n}-U\bigl(\frac{n}{\theta ks} \bigr)}{a\bigl(\frac{n}{\theta ks}\bigr)}\\
\nonumber&=& \frac{X_{n-[\theta k s],n}-U\bigl(\frac{n}{\theta k} \bigr)}{a\bigl(\frac{n}{\theta k}\bigr)} +\biggl( \frac{a\bigl(\frac{n}{\theta k}\bigr)}{a\bigl(\frac{n}{\theta ks} \bigr)}-1\biggr)\, \frac{X_{n-[\theta k s],n}-U\bigl(\frac{n}{\theta k} \bigr)}{a\bigl(\frac{n}{\theta k}\bigr)}+ \frac{U\bigl(\frac{n}{\theta k}\bigr)-U\bigl(\frac{n}{\theta ks}\bigr)}{a\bigl(\frac{n}{\theta ks}\bigr)}
\end{eqnarray}
Bearing on \eqref{RepIntOrder} combined with the uniform inequalities in Lemma \ref{LemStem}(1) and the ones for $\Pi$-varying functions provided in Proposition B.2.17 of de Haan and Ferreira (2006), we thus get for any $\varepsilon >0$,
\begin{eqnarray}\label{OpTheta}
\nonumber R_{\theta}(s) &=&-\log s +\frac{1}{\sqrt{\theta k}}\,\Bigl(\frac{W_n(s)}{s}+s^{-1/2-\varepsilon}o_p(1)\Bigr)\\
\nonumber & & \qquad \qquad \pm (s^{-\varepsilon}-1)\bigl( 1+o_p(1)\bigr)\pm \varepsilon s^{-\varepsilon}(-\log s)+\log s \pm \varepsilon s^{-\varepsilon}\\
	&=& \frac{1}{\sqrt{\theta k}}\,\frac{W_n(s)}{s}\pm (s^{-\varepsilon}-1)\bigl( 1+o_p(1)\bigr)\mp \varepsilon s^{-\varepsilon}\log s,
\end{eqnarray}
for $s \in [(\theta k)^{-1},\,1]$, all $\theta \geq 1$. Therefore, we have just seen that the distribution of deviations between high (large) sample quantiles and their theoretical counterparts is attainable with a different normalization than in \eqref{RepIntOrder}.

Before we proceed we shall require the following lemma regarding a second order condition on the auxiliary function $a$:
\begin{lem}\label{Lem2q}
Let $U\in \Pi(a)$ such that $U(\infty)= \lim_{t\rightarrow \infty}U(t)$ exists finite. Then the following limit holds with $q(t):= \int_{t}^{\infty}a(s)\,ds/s$ (defined in \eqref{qDef}),
\begin{equation*}
\limit{t}\, \frac{\frac{a(tx)}{a(t)}-1}{\frac{a(t)}{q(t)}}= -\log x, \quad x>0.
\end{equation*}
\end{lem}
\begin{pf}
The assumption that $U\in \Pi(a)$ entails
\begin{eqnarray}\label{aux1}
\nonumber    \frac{q(t)}{a(t)}\Bigl(\frac{a(tx)}{a(t)}-1\Bigr)&=& \frac{q(t)}{U(tx)-U(t)}\, \frac{U(tx)-U(t)}{a(t)}\Bigl(\frac{a(tx)}{a(t)}-1\Bigr)\\
    &=& \frac{q(t)}{U(tx)-U(t)}\,\log x\,\Bigl(\frac{a(tx)}{a(t)}-1\Bigr)\bigl(1+o(1) \bigr). \quad (t\rightarrow \infty)
\end{eqnarray}
Furthermore, according to definition \eqref{qDef} of the function $q$ and the main relation \eqref{MainRel},
\begin{equation*}
    \frac{q(t)}{U(tx)-U(t)}= \frac{\intinf{t}a(s)\, \frac{ds}{s}}{\intab{t}{tx}a(s)\, \frac{ds}{s}\bigl( 1+o(1)\bigr)}= 1+\frac{\intinf{tx}a(s)\, \frac{ds}{s}}{ \intab{t}{tx}a(s)\, \frac{ds}{s}}\bigl( 1+o(1)\bigr).
\end{equation*}
By taking the limit of the latter term when $t\rightarrow \infty$, we get from Cauchy's rule together with the fundamental theorem of integral calculus that
\begin{equation*}
    \limit{t}\frac{\intinf{tx}a(s)\, \frac{ds}{s}}{ \intab{t}{tx}a(s)\, \frac{ds}{s}}= \limit{t} \frac{-a(tx)}{a(tx)-a(t)}= -\limit{t}\Bigl( \frac{a(tx)}{a(t)}-1 \Bigr)^{-1}.
\end{equation*}
Giving heed to \eqref{aux1}, the limiting statement follows in a straightforward manner:
\begin{eqnarray*}
    \frac{q(t)}{a(t)}\biggl(\frac{a(tx)}{a(t)}-1\biggr)&=& \log x\,\Bigl(1+\frac{\intinf{tx}a(s)\, \frac{ds}{s}}{ \intab{t}{tx}a(s)\, \frac{ds}{s}}\Bigr)\Bigl(\frac{a(tx)}{a(t)}-1\Bigr)\bigl(1+o(1) \bigr)\\
        &=& -\log x  +\log x \Bigl(\frac{a(tx)}{a(t)}-1\Bigr)\bigl(1+o(1) \bigr). \quad (t\rightarrow \infty)
\end{eqnarray*}
\end{pf}

\begin{prop}\label{PropConsQ}
Let $X_1,\,X_2,\ldots$ be i.i.d. random variables with tail quantile function $U$ satisfying condition \eqref{MainRel}. Suppose  $k=k(n)$ is a sequence of positive integers such that $k(n)\rightarrow \infty$, $k(n)/n\rightarrow 0$, as $n\rightarrow \infty$. Then $\hat{q}(n/k)$ is a consistent estimator for $q(n/k)$ in the sense that the following convergence in probability holds,
\begin{equation*}
\frac{\hat{q}\bigl(\ndivk \bigr)}{q\bigl(\ndivk \bigr)}\conv{p}1.
\end{equation*}
\end{prop}
\begin{pf}
We begin by noting that
\begin{eqnarray}
\nonumber \frac{\hat{q}\bigl(\ndivk \bigr)}{q\bigl(\ndivk \bigr)}& =& -\frac{1}{\log 2} \intab{0}{1}\frac{\widehat{U}\bigl(\frac{n}{2ks} \bigr)-\widehat{U}\bigl(\frac{n}{ks} \bigr)}{q\bigl( \ndivk\bigr)} \, \frac{ds}{s}\\
\label{q1overq}  & = & -\frac{1}{\log 2}\,  \biggl\{\intab{\frac{1}{2k}}{1}\frac{X_{n-[2ks],n}-U\bigl(\frac{n}{2ks} \bigr)}{q\bigl(\ndivk \bigr)}\, \frac{ds}{s}-\intab{\frac{1}{k}}{1}\frac{X_{n-[ks],n}-U\bigl(\frac{n}{ks} \bigr)}{q\bigl( \ndivk\bigr)} \, \frac{ds}{s}\\
\label{q1overqNext}& & \mbox{\hspace{2.5cm}} -\intab{\frac{1}{2k}}{\frac{1}{k}} \frac{X_{n,n} -U\bigl(\frac{n}{ks} \bigr)}{q\bigl( \ndivk\bigr)} \,\frac{ds}{s}+ \intab{\frac{1}{2k}}{1}\frac{U\bigl(\frac{n}{2ks}\bigr)-U\bigl(\frac{n}{ks} \bigr)}{q\bigl( \ndivk\bigr)}\, \frac{ds}{s}\biggr\}.
\end{eqnarray}
The two integral terms in \eqref{q1overq} shall be handled jointly through the consideration of $R_{2}(s)$ (see Eq. \eqref{RTheta} with $\theta =2$) in the one integral below:
\begin{equation}\label{2in1}
\intab{\frac{1}{2k}}{1}\frac{X_{n-[2ks],n}-U\bigl(\frac{n}{2ks} \bigr)}{q\bigl(\ndivk \bigr)}\, \frac{ds}{s}-\intab{\frac{1}{k}}{1}\tfrac{X_{n-[ks],n}-U\bigl(\frac{n}{ks} \bigr)}{q\bigl( \ndivk\bigr)} \, \frac{ds}{s}=\intab{\frac{1}{2}}{1} \tfrac{X_{n-[2k s],n}-U\bigl(\frac{n}{2 ks} \bigr)}{q\bigl(\ndivk \bigr)}\, \frac{ds}{s}=:I_1(k,n)
\end{equation}
whence
\begin{equation}\label{2in1Exp}
I_1(k,n)= \frac{a\bigl(\ndivk\bigr)}{q\bigl( \ndivk\bigr)} \Bigl\{ \intab{\frac{1}{2}}{1}R_{2}(s)\, \frac{ds}{s} +  \intab{\frac{1}{2}}{1}\Bigl( \frac{a\bigl(\frac{n}{2 ks}\bigr)}{a\bigl(\ndivk\bigr)}-1\Bigr) R_{2}(s)\, \frac{ds}{s}\Bigr\}.
\end{equation}
Now, Lemma \ref{Lem2q} ascertains
\begin{eqnarray}\label{ubound}
\nonumber I_1(k,n)&=& \intab{\frac{1}{2}}{1} \frac{X_{n-[2k s],n}-U\bigl(\frac{n}{2 ks} \bigr)}{q\bigl(\ndivk \bigr)}\, \frac{ds}{s}\\
\nonumber &\leq &\frac{a\bigl(\ndivk\bigr)}{q\bigl(\ndivk \bigr)}\Bigl|\intab{\frac{1}{2}}{1} R_{2}(s)\, \frac{ds}{s}\Bigr|+\Bigl(\frac{a\bigl(\ndivk\bigr)}{q\bigl(\ndivk\bigr)}\Bigr)^2 \Bigl|\intab{\frac{1}{2}}{1} R_{2}(s) \log(2 s)\, \frac{ds}{s}\Bigr|\\
&\leq &\frac{a\bigl(\ndivk\bigr)}{q\bigl(\ndivk \bigr)}\biggl(1+\frac{a\bigl(\ndivk\bigr)}{q\bigl(\ndivk\bigr)}\log 2\biggr)\Bigl|\intab{\frac{1}{2}}{1} R_{2}(s)\, \frac{ds}{s}\Bigr|,
\end{eqnarray}
with high probability, for sufficiently large $n$. We can provide a similar lower bound.

\noindent Owing to \eqref{OpTheta}, the following holds w.r.t. the integral featuring in the upper bound \eqref{ubound}, for any positive $\varepsilon$,
\begin{equation*}
\Bigl|\intab{\frac{1}{2}}{1}R_{2}(s)\, \frac{ds}{s}\Bigr|\leq  \Bigl|\frac{1}{\sqrt{2 k}}\intab{\frac{1}{2}}{1} s^{-2}\,W_n(s) \, ds\Bigr| + \intab{\frac{1}{2}}{1}\bigl(s^{-\varepsilon}-1\bigr)\, \frac{ds}{s}\bigl(1+o_p(1)\bigr) -\varepsilon \intab{\frac{1}{2}}{1}s^{-\varepsilon} \log s
\, \frac{ds}{s}.
\end{equation*}
Since $\varepsilon>0$ is arbitrary, then
\begin{equation*}
	0< \intab{\frac{1}{2}}{1} \bigl(s^{-\varepsilon}-1\bigr)\, \frac{ds}{s}= \frac{2^{\varepsilon}-1}{\varepsilon}-\log 2\, \arrowd{\varepsilon \downarrow 0} 0,
	\end{equation*}
meaning that
\begin{equation*}
\intab{\frac{1}{2}}{1} \bigl(s^{-\varepsilon}-1\bigr)\, \frac{ds}{s}
\end{equation*}
can be discarded. A similar line of reasoning applies to
\begin{equation*}
\varepsilon \intab{\frac{1}{2}}{1}s^{-\varepsilon} \log \bigl(\frac{1}{s}\bigr)
\, \frac{ds}{s}=2^{\varepsilon}\log 2-\frac{2^{\varepsilon}-1}{\varepsilon}\, \arrowd{\varepsilon \downarrow 0} 0,
\end{equation*}
thus also discarded.

\noindent We now recall that $k=k(n)$ is a sequence of positive integers tending to infinity as $n\rightarrow \infty$. Let us define
\begin{equation*}
   Y_n:=\frac{1}{\sqrt{2 k}}\intab{\frac{1}{2}}{1} W_n(s)\, \frac{ds}{s^{2}},
\end{equation*}
which regards a sequence of normal random variables with zero mean and variance equal to
\begin{equation*}
Var(Y_n)=\frac{1-\log 2}{k} \; \arrowf{n} 0.
\end{equation*}
The latter means that the sequence of random variables $\{Y_n\}_{n\geq 0}$ is a sequence of degenerate random variables, eventually, and
the two integrals in \eqref{q1overq}  (unifyed in \eqref{2in1}; see also Eq. \eqref{2in1Exp} in terms of  $R_2(s)$) vanish  with probability tending to one as $n\rightarrow \infty$. In this respect we note that $a(n/k)/q(n/k)=o(1)$, which entails in fact that
\begin{equation*}
I_1(k,n)=o_p(1)\,\,\left(=o_p(\tfrac{a(n/k)}{q(n/k)})=O_p(\tfrac{a(n/k)}{\sqrt{k}q(n/k)})\right).
\end{equation*}

\noindent The rest of the proof pertains to the terms in \eqref{q1overqNext}. Regarding the first integral in \eqref{q1overqNext}, we note that
\begin{eqnarray*}
\nonumber I_2(k,n)&:=&\intab{\frac{1}{2k}}{\frac{1}{k}}\frac{X_{n,n}-U\bigl(\frac{n}{ks}\bigr)}{q\bigl(\ndivk \bigr)}\, \frac{ds}{s}\\
\nonumber &=& \intab{1/2}{1} \frac{X_{n,n}-U\bigl(\frac{n}{s}\bigr)}{q\bigl(\ndivk \bigr)}\, \frac{ds}{s}\\
\nonumber & \id& \frac{a(n)}{q\bigl(\ndivk\bigr)} \biggl\{\frac{U\bigl(\frac{1}{U_{1,n}}\bigr)-U(n)}{a(n)}\log 2  -\intab{\frac{1}{2}}{1}\frac{U\bigl(\frac{n}{s}\bigr)-U(n)}{a(n)}\, \frac{ds}{s}\biggr\}\\
\nonumber &= & \frac{a(n)}{q\bigl(\ndivk\bigr)} \Bigl\{ -\log (n\,U_{1,n})\log 2 + \intab{\frac{1}{2}}{1} \log s\,\frac{ds}{s} \Bigr\}\bigl(1+o_p(1)\bigr)\\
&=& \frac{a(n)}{q\bigl(\ndivk\bigr)}\log 2 \Bigl( -\log (n\,U_{1,n}) - \frac{1}{2}\log 2\Bigr)\bigl(1+o_p(1)\bigr).
\end{eqnarray*}

\noindent Now, the probability integral transformation yields the following equality in distribution for the random term above:
\begin{equation}\label{Max}
-\log (n\, U_{1,n})\,\id \,E_{n,n}-\log n,
\end{equation}
where $E_{n,n}$ is the maximum of $n$ i.i.d. standard exponential random variables. Hence, the random variable  \eqref{Max} converges in distribution to a Gumbel random variable with distribution function given by $\exp\{-e^{-x}\}$, $x\in \real$. Moreover, $a(n)/q(n/k)\rightarrow 0$, as $n\rightarrow \infty$, because $a(n/k)/q(n/k)=o(1)$ (see Lemma \ref{LemStem}(2)), where the auxiliary positive function $a$ satisfies  $a(t)\rightarrow 0$, as $t\rightarrow \infty$, by assumption. Therefore,
\begin{equation}\label{I2}
I_2(k,n)=o_p(1)\,\,\left(=O_p(\tfrac{a(n)}{q(n/k)}) \right).
\end{equation}
\noindent In order to finally attain consistency of $\hat{q}(n/k)$ let us consider the last integral in \eqref{q1overqNext}, which we will show it is bounded. On the one hand, for the upper bound,
\begin{equation}\label{UpBInt}
\intab{\frac{1}{2 k}}{1}\frac{U\bigl(\frac{n}{ks}\bigr)-U\bigl(\frac{n}{2ks} \bigr)}{q\bigl(\ndivk\bigr)}\, \frac{ds}{s}\leq  \intunit\frac{U\bigl(\frac{n}{ks}\bigr)-U\bigl(\frac{n}{2ks} \bigr)}{q\bigl(\ndivk\bigr)}\, \frac{ds}{s}
\end{equation}
and on the other hand, for the lower bound,
\begin{eqnarray*}
	& & \intab{\frac{1}{2 k}}{1}\frac{U\bigl(\frac{n}{ks}\bigr)-U\bigl(\frac{n}{2ks} \bigr)}{q\bigl(\ndivk\bigr)}\, \frac{ds}{s}\\
         &=& \intab{0}{1-\frac{1}{2k}}\frac{U\bigl(\frac{n}{ks+1/2}\bigr)-U\bigl(\frac{n}{2ks+1}\bigr)}{q\bigl(\frac{n}{k}\bigr)}\, \frac{ds}{s+\frac{1}{2k}} \\
	&\geq &\intunit\frac{  U\Bigl(\frac{n}{k\bigl(s+\frac{1}{2k}\bigr)}\Bigr)- U\Bigl(\frac{n}{2k\bigl(s+\frac{1}{2k}\bigr)}\Bigr)}{q\bigl(\frac{n}{k}\bigr)}\,\frac{ds}{s+\frac{1}{2k}}-2\intab{1-\frac{1}{2k}}{1}\frac{U\bigl(\frac{n}{ks+1/2}\bigr)-U\bigl(\frac{n}{2ks+1}\bigr)}{q\bigl(\frac{n}{k}\bigr)}\, \frac{ds}{s+\frac{1}{2k}}.
\end{eqnarray*}
Making $t=n/k$ run on the real line towards infinity, then the $\Pi-$variation in relation \eqref{PiInt} is rephrased as
\begin{equation}\label{RephPiInt}
\limit{t} \frac{\intunit U\bigl(\frac{tx}{s}\bigr)\,\frac{ds}{s}-\intunit U\bigl(\frac{t}{s}\bigr)\,\frac{ds}{s}}{q(t)}= \log x, \quad x>0,
\end{equation}
which clearly entails the following limit for the upper bound in \eqref{UpBInt}:
\begin{equation*}
   \frac{ \intunit U\bigl(\frac{n}{ks}\bigr)\, \frac{ds}{s}- \intunit U\bigl(\frac{n}{2ks} \bigr)\, \frac{ds}{s}}{q\bigl(\ndivk\bigr)}=- \frac{ \intunit U\bigl(\frac{n}{2ks}\bigr)\, \frac{ds}{s}- \intunit U\bigl(\frac{n}{ks} \bigr)\, \frac{ds}{s}}{q\bigl(\ndivk\bigr)}\arrowf{n} \log 2.
\end{equation*}
Regarding the lower bound,
\begin{eqnarray}
\label{LowBInt1}	\intab{\frac{1}{2 k}}{1}\frac{U\bigl(\frac{n}{ks}\bigr)-U\bigl(\frac{n}{2ks} \bigr)}{q\bigl(\ndivk\bigr)}\, \frac{ds}{s}&\geq &\intunit \frac{ U\Bigl(\frac{n}{k\bigl(s+\frac{1}{2k}\bigr)}\Bigr)- U\Bigl(\frac{n}{2k\bigl(s+\frac{1}{2k}\bigr)}\Bigr)}{q\bigl(\frac{n}{k}\bigr)}\,\frac{ds}{s+\frac{1}{2k}}\\
\label{LowBInt2}	& & \mbox{\hspace{0.5cm}}-2\frac{q\bigl(\frac{n}{2k}\bigr)}{q\bigl(\frac{n}{k}\bigr)}\intab{1-\frac{1}{2k}}{1}\frac{U\bigl(\frac{n}{ks+1/2}\bigr)-U\bigl(\frac{n}{2ks+1}\bigr)}{q\bigl(\frac{n}{2k}\bigr)}\, \frac{ds}{s+\frac{1}{2k}},
\end{eqnarray}
we note that for every $\varepsilon>0$, there exists $n_0 \in \field{N}$ such that for $n\geq n_0$,
\begin{equation}
	\Bigl| \frac{1}{s+1/(2k)}-\frac{1}{s}\Bigr|<\varepsilon.
\end{equation}
Whence, we have in turn the following inequality with respect to \eqref{LowBInt1}:
\begin{equation*}
	\intunit \frac{ U\Bigl(\frac{n}{k\bigl(s+\frac{1}{2k}\bigr)}\Bigr)- U\Bigl(\frac{n}{2k\bigl(s+\frac{1}{2k}\bigr)}\Bigr)}{q\bigl(\frac{n}{k}\bigr)}\,\frac{ds}{s+\frac{1}{2k}}> \intunit\frac{  U\Bigl(\frac{n}{k\bigl(s+\frac{1}{2k}\bigr)}\Bigr)- U\Bigl(\frac{n}{2k\bigl(s+\frac{1}{2k}\bigr)}\Bigr)}{q\bigl(\frac{n}{k}\bigr)}\Bigl(\frac{1}{s}-\varepsilon\Bigr)\, ds.
\end{equation*}
For the first part of the right-hand side of the above we use again condition \eqref{RephPiInt}, while the second part is dealt with Theorem B.2.19 of \citet{deHaanFerreira:06} involving the fact that $U\in \Pi(a)$:
\begin{eqnarray*}
 & &\intunit\frac{  U\Bigl(\frac{n}{k\bigl(s+\frac{1}{2k}\bigr)}\Bigr)- U\Bigl(\frac{n}{2k\bigl(s+\frac{1}{2k}\bigr)}\Bigr)}{q\bigl(\frac{n}{k}\bigr)}\,\frac{ds}{s}+\varepsilon\frac{a\bigl(\frac{n}{k}\bigr)}{q\bigl(\frac{n}{k}\bigr)}\intunit\frac{ U\Bigl(\frac{n}{2k\bigl(s+\frac{1}{2k}\bigr)}\Bigr)- U\Bigl(\frac{n}{k\bigl(s+\frac{1}{2k}\bigr)}\Bigr)}{a\bigl(\frac{n}{k}\bigr)}\, ds\\
 &=& \log 2 \bigl(1+o(1)\bigr)-\varepsilon \frac{a\bigl(\frac{n}{k}\bigr)}{q\bigl(\frac{n}{k}\bigr)}\log 2 \bigl(1+o(1)\bigr)\\
 & \arrowf{n}& \log 2.
\end{eqnarray*}
For the latter, we recall that $a(n/k)=o\bigl(q(n/k)\bigr)$.

\noindent Now we write $\delta=1/(2k)>0$ everywhere   in \eqref{LowBInt2}.  Furthermore, we assume that there exists $n_0\in \field{N}$ such that, for $n\geq n_0$, the term  $n\delta$ is   large enough and the integral in \eqref{LowBInt2} can rephrased as
\begin{equation}\label{Intdelta}
	 I^*_{\delta}:=\frac{\intab{1-\delta}{1}\Bigl(U\bigl(\frac{2}{s+\delta}n\delta\bigr)-U\bigl(\frac{1}{s+\delta}n\delta\bigr)\Bigr)\, \frac{ds}{s+\delta}}{\intab{n\delta}{1}a(s)\, \frac{ds}{s}}.
\end{equation}
We note that, for every fixed $\delta >0$, we have that from the $\Pi$-variation of $U$ that the following holds for the numerator of $I^*_{\delta}$ properly rescaled by $a(n\delta)$ \citep[cf. Theorem B.2.19 in][]{deHaanFerreira:06}:
\begin{equation*}
	  \frac{\intab{1-\delta}{1}\Bigl(U\bigl(\frac{2}{s+\delta}n\delta\bigr)-U\bigl(\frac{1}{s+\delta}n\delta\bigr)\Bigr)\, \frac{ds}{s+\delta}}{a(n\delta)}\arrowf{n}\intab{1-\delta}{1}\log 2\,\frac{ds}{s+\delta}= \log(1+\delta)\log2.
\end{equation*}
For arbitrary small $\delta$, the latter approaches zero. Predicated on the above, we apply Cauchy's rule to obtain $\lim_{\delta\rightarrow 0}I^*_{\delta}$ (we recall that $\delta\rightarrow 0$ implies $n\rightarrow \infty$). Towards this end, we apply Eq. (2.11) of \citet{Chiang:00} upon the numerator of $I^*_{\delta}$, whence
\begin{eqnarray*}
	\lim_{\delta\rightarrow 0} I^*_{\delta}&=&\lim_{\delta\rightarrow 0} \frac{\intab{1-\delta}{1}\Bigl(U'\big(\frac{2}{s+\delta}n\delta\bigr)\frac{2s}{(s+\delta)^3}-U'\big(\frac{1}{s+\delta}n\delta\bigr)\frac{s}{(s+\delta)^3}\Bigr)\,ds}{-\frac{a(n\delta)}{n\delta}}\\
	& &\mbox{\hspace{1.0cm}} +\lim_{\delta\rightarrow 0}\Bigl\{\delta\intab{1-\delta}{1}\frac{U\bigl(\frac{2n\delta}{s+\delta}\bigr)-U\bigl(\frac{n\delta}{s+\delta}\bigr)}{a(n\delta)}\,\frac{ds}{(s+\delta)^2}- \delta\frac{U(2n\delta)-U(n\delta)}{a(n\delta)} \Bigr\}.
\end{eqnarray*}
Since $U'(t)=a(t)/t$ then the limit becomes equal to the the limit of
\begin{equation*}
-\intab{1-\delta}{1} \Bigl(\frac{a\bigl(\frac{2n\delta}{s+\delta}\bigr)}{a(n\delta)}-\frac{a\bigl(\frac{n\delta}{s+\delta}\bigr)}{a(n\delta)}\Bigr)\,
\frac{s\,ds}{(s+\delta)^2}+\delta\Bigl(\intab{1-\delta}{1}\frac{U\bigl(\frac{2n\delta}{s+\delta}\bigr)-U\bigl(\frac{n\delta}{s+\delta}\bigr)}{a(n\delta)}\,\frac{ds}{(s+\delta)^2}- \frac{U(2n\delta)-U(n\delta)}{a(n\delta)}\Bigr).
\end{equation*}
We can now take any arbitrary small $\delta $ (making $n\rightarrow \infty$) in order to apply the uniform convergence of $a\in RV_0$ and $U\in\Pi(a)$ so that the above integrals are ensured finite and then equal to zero by definition. Hence, all the terms are negligible as $\delta$ converges to zero meaning that $\lim_{\delta \rightarrow 0}I^*_{\delta}$ becomes null.
Therefore,
\begin{equation*}
\intab{\frac{1}{2 k}}{1}\frac{U\bigl(\frac{n}{2ks}\bigr)-U\bigl(\frac{n}{ks} \bigr)}{q\bigl(\ndivk\bigr)}\, \frac{ds}{s} \arrowf{n} -\log 2.
\end{equation*}
and the precise result for consistency of $\hat{q}(n/k)$ thus follows by noting that $q(n/k)\sim q\bigl(n/(2k)\bigr)$.
\end{pf}

In view of \eqref{EstEndPoint}, we have the following alternative formulation aimed at establishing consistency of the proposed estimator for the right endpoint.

\begin{thm}\label{TheoConsxF}
Let $X_1,\,X_2,\ldots$ be i.i.d. random variables with tail quantile function $U$ satisfying condition \eqref{MainRel}. Suppose  $k=k(n)$ is a sequence of positive integers such that $k(n)\rightarrow \infty$, $k(n)/n\rightarrow 0$, as $n\rightarrow \infty$. Then $\hat{x}^{F}:= X_{n-k,n}+ \hat{q}(n/k)$ is a consistent estimator for $x^{F}<\infty$, i.e.
\begin{equation*}
\hat{x}^F\conv{p}x^F.
\end{equation*}
\end{thm}
\begin{pf}
It will suffice to note there are three main contributing components for $x^{F}-\hat{x}^F$. Specifically,
\begin{eqnarray*}
x^F-\hat{x}^F &\id& U(\infty)-U\bigl(\frac{1}{U_{k+1,n}}\bigr)-\hat{q}\bigl(\ndivk\bigr)\\
		&=& \Bigl(U(\infty)- U\bigl(\ndivk\bigr)-q\bigl(\ndivk\bigr)\Bigr)-\Bigl(U\bigl(\frac{1}{U_{k+1,n}}\bigr)-U\bigl(\ndivk\bigr)\Bigr)-q\bigl(\ndivk\bigr)\Bigl(\frac{\hat{q}\bigl(\ndivk\bigr)}{q\bigl(\ndivk\bigr)}-1\Bigr)\\
		&=& I-II-III,		
\end{eqnarray*}
where:
\begin{equation*}
 I:=U(\infty)- U\bigl(\ndivk\bigr)-q\bigl(\ndivk\bigr)=o\Bigl(a\bigl(\ndivk\bigr)\Bigr),
\end{equation*}
which follows directly from relation \eqref{MainRel};
\begin{equation*}
II:= U\bigl(\frac{1}{U_{k+1,n}}\bigr)-U\bigl(\ndivk\bigr)=o_p\Bigl(a\bigl(\ndivk\bigr)\Bigr)
\end{equation*}
because $U\in\Pi(a)$ while Smirnov's Lemma ensures $k/(n\,U_{k+1,n})\conv{P} 1$ \cite[see Lemma 2.2.3 in][]{deHaanFerreira:06};
\begin{equation*}
III:= q\bigl(\ndivk\bigr)\Bigl(\frac{\hat{q}\bigl(\ndivk\bigr)}{q\bigl(\ndivk\bigr)}-1\Bigr)=o_p(1)
\end{equation*}
which is verified by Proposition \ref{PropConsQ} and the fact that relation \eqref{MainRel} implies $q(n/k)=o(1)$.

\end{pf}

\vspace{0.2cm}%

The asymptotic distribution of $\hat{q}(n/k)$ is predicated on a suitable second order refinement of \eqref{PiofU}: suppose there exist functions $a$, positive and $A$, positive or negative, both tending to zero as $t\rightarrow \infty$, such that
\begin{equation}\label{2ndPiofU}
\limit{t}\frac{\frac{U(tx)-U(t)}{a(t)}-\log x}{A(t)}= \frac{1}{2}(\log x)^2,
\end{equation}
for all $x>0$.

\begin{rem}
The second order condition above follows directly from Theorem B.3.6, Remark B.3.7 and Corollary 2.3.5 of \cite{deHaanFerreira:06} because the former states that, in our setup of $\gamma=0$ and $x^F<\infty$, the only case allowed is the case of the second order parameter $\rho$ equal to zero. Like the function $a$, the second order auxiliary function $A$ converges to zero, not changing sign for $t$ near infinity, and $|A|$ is slowly varying, i.e. $A(tx)/A(t) \rightarrow 1, t\rightarrow \infty$ (notation: $|A|\in RV_0$).
\end{rem}

Furthermore, Theorem 2.3.6 of  \cite{deHaanFerreira:06} ascertains the existence of functions $a_0$ and $A_0$ satisfying, as $t\rightarrow \infty$, $A_0(t) \sim A(t)$ and $a_0(t)/a(t)-1=o\bigl(A(t)\bigr)$, with the property that for any $\varepsilon >0$, there exists $t_0=t_0(\varepsilon)$ such that for all $t, \, tx \geq t_0$,
\begin{equation}\label{2ndUnifBounds}
\biggl| \frac{\frac{U(tx)-U(t)}{a_0(t)}-\log x}{A_0(t)}-\frac{1}{2}(\log x)^2\biggr|\leq \varepsilon \max(x^{\varepsilon}, x^{-\varepsilon})
\end{equation}
and
\begin{equation}\label{AuxUnifBounds}
\biggl| \frac{\frac{a_0(tx)}{a_0(t)}-1}{A_0(t)}-\log x\biggr|\leq \varepsilon \max(x^{\varepsilon}, x^{-\varepsilon}).
\end{equation}

\begin{rem}\label{Rem2}
We note that relation \eqref{AuxUnifBounds} combined with Lemma \ref{Lem2q} ascertains that $-a_0(t)/q(t)=c A_0(t)$, with $c\neq 0$ because $\rho=\gamma=0$ \citep[cf. Eq. (B.3.4) and Remark B.3.5 in][]{deHaanFerreira:06}. Henceforth we may assume that the function $q$ is conveniently redefined so that $-a/q \sim A$.
\end{rem}

\begin{ex}\label{ExNegFrech}
The Negative Fr\'echet model with parameter $\beta>0$, \ie, with distribution function $F(x)= 1-\exp\{-(x^F-x)^{-\beta}\}$, $x\geq x^F$, $\beta>0$. The associated tail quantile function $U$ is given by  $U(t)= x^F-(\log t)^{-1/\beta}$, $t\geq 1$. Then $U\in\Pi(a_0)$ with $a_0(t)= (1/\beta) (\log t)^{-1/\beta-1}\rightarrow 0$, as $t\rightarrow \infty$. Therefore, the auxiliary function $q$ defined in \eqref{qDef} becomes $q(t)= (\log t)^{-1/\beta}$, $\beta>0$. Now, by straightforward  calculations we see that $A_0(t)=-(1+1/\beta)(\log t)^{-1}$, which implies that $-a_0(t)/q(t)=A_0(t)/(1+\beta)$, for $t$ near infinity.
\end{ex}

Theorem 2.4.2 of \cite{deHaanFerreira:06} allows to gain insight about the distributional representation displayed in \eqref{RepIntOrder}. Specifically, if the tail quantile function satisfies the second order condition \eqref{2ndPiofU} then, for each $\varepsilon >0$,
\begin{equation}\label{TheodH&F}
\suprem{\frac{1}{\theta k}\leq s\leq 1} s^{1/2+\varepsilon} \biggl| \sqrt{\theta k}\biggl(\frac{X_{n-[\theta k s],n}-U\bigl(\frac{n}{\theta k}\bigr)}{a_0\bigl(\frac{n}{\theta k}\bigr)}+ \log s\biggr)- \frac{W_n(s)}{s}-\sqrt{\theta k}\,A_0\bigl( \frac{n}{\theta k}\bigr)\frac{1}{2}(\log s)^2\biggr|\conv{p}0,
\end{equation}
provided $k=k(n) \rightarrow \infty$, $k/n=o(n)$ and $\sqrt{k}A_0(n/k)=O(1)$.

Therefore, the asymptotic distribution of $\hat{q}(n/k)$ will appear intertwined with the proof of consistency in Proposition \ref{PropConsQ} via $R_{\theta}(s)$, (defined in \eqref{RTheta} for $s\in [(\theta k)^{-1},1]$, see also \eqref{OpTheta}), albeit under the second order grasp provided above. The next Proposition accounts for this \cite[cf. (2.4.7) of][]{deHaanFerreira:06}.
\begin{prop}\label{PropRTheta}
Suppose the second order condition \eqref{2ndPiofU} holds. Let $k=k(n) \rightarrow \infty$, $k/n=o(n)$ and $\sqrt{k}A(n/k)\rightarrow \lambda \in \real$, as $n\rightarrow \infty$. Then, for $\theta \geq 1$ and for each $\varepsilon>0$ sufficiently small,
\begin{equation*}
\suprem{\frac{1}{\theta k}\leq s\leq 1} s^{1/2+\varepsilon} \biggl|\sqrt{\theta k}\,\frac{X_{n-[\theta k s],n}-U\bigl(\frac{n}{\theta ks} \bigr)}{a\bigl(\frac{n}{\theta ks}\bigr)} -\frac{ W_n(s)}{s}\biggr|=o_p(1).
\end{equation*}

\end{prop}
\begin{pf}
Similarly to the equality right after \eqref{RTheta}, we have that
\begin{equation*}
R_{\theta}(s):=\,\frac{X_{n-[\theta k s],n}-U\bigl(\frac{n}{\theta ks} \bigr)}{a\bigl(\frac{n}{\theta ks}\bigr)}\\
= \frac{a_0\bigl(\frac{n}{\theta k }\bigr)}{a\bigl(\frac{n}{\theta k s} \bigr)} \, \Biggl\{\frac{X_{n-[\theta k s],n}-U\bigl(\frac{n}{\theta k} \bigr)}{a_0\bigl(\frac{n}{\theta k}\bigr)} - \frac{U\bigl(\frac{n}{\theta k s}\bigr)-U\bigl(\frac{n}{\theta k}\bigr)}{a_0\bigl(\frac{n}{\theta k}\bigr)}\Biggr\}.
\end{equation*}
Noting that
\begin{equation*}
\frac{a_0(t)}{a\bigl(\frac{t}{s}\bigr)}= \frac{a_0(t)}{a(t)}\,\frac{a\bigl(t\bigr)}{a\bigl(\frac{t}{s}\bigr)},
\end{equation*}
for all $s>0$, then Lemma \ref{Lem2q} combined with Remark \ref{Rem2} yields the expansion
\begin{equation}\label{AuxExp}
\frac{a_0(t)}{a\bigl(\frac{t}{s}\bigr)}= \frac{a_0(t)}{a(t)}\, \Biggl(1-\frac{a( t)}{q(t)}\log s+ o\Bigl(\frac{a( t)}{q(t)}\Bigr)\Biggr)= \frac{a_0(t)}{a(t)}\,\Bigl( 1+A(t)\log s+o\bigl(A(t)\bigr)\Bigr),
\end{equation}
for all $s>0$. In this respect, we also note that $|A|\in RV_0$ and $a_0(t)/a(t)= 1+o\bigl(A(t)\bigr)$.

\noindent Having set $ 1/(\theta k) \leq s\leq 1$, we thus have from \eqref{TheodH&F}, the uniform bounds in \eqref{2ndUnifBounds} and the second equality in \eqref{AuxExp}, that
\begin{eqnarray*}
\sqrt{\theta k}\,R_{\theta}(s)&=& \frac{W_n(s)}{s}+A\bigl(\frac{n}{\theta k}\bigr)\frac{\log s}{s} W_n(s) \mp \varepsilon s^{-\varepsilon}\sqrt{\theta k}\, A\big(\frac{n}{\theta k}\bigr)\pm \varepsilon s^{-\varepsilon}\log s \,\sqrt{\theta k}\, A^2\big(\frac{n}{\theta k}\bigr)\\
	& & \mbox{\hspace{5.0cm}} +o_p(s^{-\frac{1}{2}-\varepsilon})+ \,o_p\Bigl(s^{-\frac{1}{2}-\varepsilon}\log s\,A\bigl(\frac{n}{\theta k}\bigr) \Bigr),
\end{eqnarray*}
uniformly in $s$.
Hence, the assumption that $\sqrt{k}A(n/k)=O(1)$ entails that $\log (1/s) A\bigl(n/(\theta k)\bigr)\rightarrow 0$,  whereas $\varepsilon s^{-\varepsilon} \sqrt{\theta k}A\bigl(n/(\theta k)\bigr)$ virtually becomes $o(s^{-1/2-\varepsilon})$ for each $\varepsilon>0$ arbitrarily small and uniformly in $s\in[(\theta k)^{-1},1]$. The $o_p$-terms are uniform in $s \in [1/(\theta k),1]$. Hence the following representation for $\sqrt{\theta k}\,R_{\theta}(s)$, valid for $\varepsilon \in (0,1)$,
\begin{equation*}
\sqrt{\theta k}\,R_{\theta}(s)= \frac{W_n(s)}{s}+o_p(s^{-1/2-\varepsilon}).
\end{equation*}
\end{pf}

\vspace{0.2cm}
\begin{thm}\label{TheoAN}
Assume the second order condition \eqref{2ndPiofU} holds. Suppose  $k=k(n)$ is such that, as $n\rightarrow \infty$, $k(n)\rightarrow \infty$, $k(n)/n\rightarrow 0$, $a(n)/a(n/k)\rightarrow 1$ and $\sqrt{k}\,A(n/k)=O(1)$. Assume furthermore that
\begin{equation}\label{Restrict}
\limit{n}\, \frac{1}{A(n/k)}\biggl( \intab{\frac{1}{2k}}{1}\frac{U\bigl(\frac{n}{ks}\bigr)-U\bigl(\frac{n}{2ks}\bigr)}{q\bigl(\ndivk\bigr)}\, \frac{ds}{s}-\log 2\biggr)=\lambda\in \real.
\end{equation}
Then
\begin{equation}\label{Limit}
\frac{q\bigl(\ndivk\bigr)}{a\bigl(\ndivk\bigr)}\Bigl( \frac{\hat{q}\bigl(\ndivk\bigr)}{q\bigl(\ndivk\bigr)}-1\Bigr)\conv{d} \Lambda -\frac{\log 2}{2}-\frac{\lambda}{\log 2},
\end{equation}
where $\Lambda$ is a Gumbel random variable with distribution function $\exp\{-e^{-x}\}$, all $x\in \real$.
\end{thm}

Before giving a proof, we note that the assumption \eqref{Restrict} of the theorem regards a second order refinement of \eqref{PiInt}, more concretely:
 \begin{equation}\label{2ndExtRVIntU}
	\limit{t} \frac{\frac{\intinf{tx}U(s)\,\frac{ds}{s}-\intinf{t}U(s)\,\frac{ds}{s}}{q(t)}- \log x}{Q(t)}= \frac{1}{2}(\log x)^2,
\end{equation}
taken in the point $x=2$ for large enough $t=n/k$. Hence, the assumption \eqref{Restrict} has been tailored via the usual second order setup (see also Eq. \eqref{2ndPiofU}) provided by the theory of extended regular variation, with $Q(t)=O(A(t))$. We refer to Appendix B of \citet{deHaanFerreira:06} for a good catalog on results concerning theory of extended regular variation.

The assumption on that $a(n/k)/a(n) \rightarrow 1$, as $n \rightarrow \infty$ is, however, a bit more restrictive in terms of screening for an adequate value $k$ which will determine the number of top order statistics to base our inference from. For example, if we assume the Negative Fr\'echet for the underlying distribution function (see Example \ref{ExNegFrech}) and $k_n=n^{p}$, $p\in (0,1)$, then
\begin{equation*}
	\frac{a(n)}{a(n/k_n)}= \Bigl(1-\frac{\log k_n}{\log n}\Bigr)^{1/\beta+1}= (1-p)^{1/\beta+1},
\end{equation*}
which is approximately $1$ if and only if $p$ approaches zero. A more appropriate choice regards intermediate sequences at a slower rate such as $k_n= (\log n)^r$, $r\in (0,2]$. Bearing this choice in mind, we have that
\begin{equation*}
	\frac{a(n)}{a(n/k_n)}= \Bigl(1-\frac{\log k_n}{\log n}\Bigr)^{1/\beta+1}=  \Bigl(1-\frac{r}{\log n}+\frac{\log \log n}{\log n}\Bigr)^{1/\beta+1} \arrowf{n} 1.
\end{equation*}
The upper bound $r\leq 2$ is imposed in order to comply with the assumption $\sqrt{k_n}A(n/k_n)=O(1)$.

Given the slow variation feature of all the functions involved in the characterizations of the present subclass of distributions in the Gumbel domain with finite right endpoint, we believe that the latter choice for $k=k_n$ is a feasible one for most models satisfying \eqref{MainRel}, meaning that we require intermediate values $k_n$ such that
 $\log (k_n)=o(\log n)$.
Altogether, we are excluding  Nevertheless, we can bring forward the fact that a miss-specification of $k_n$ in the sense that $a(n/k_n)/a(n)$ converges to a constant different than $1$, has a direct impact on the asymptotic variance of the normalized relative error presented in Theorem \ref{TheoAN} rather than upon the asymptotic bias. This can be clearly seen in the proof we present below.

\bigskip

\begin{pfofThm}{\bf \ref{TheoAN}:}
Similarly as in \eqref{2in1}, we have that
\begin{eqnarray}
\nonumber & &\frac{q\bigl(\ndivk\bigr)}{a\bigl(\ndivk\bigr)}\Bigl( \frac{\hat{q}\bigl(\ndivk\bigr)}{q\bigl(\ndivk\bigr)}-1\Bigr) \\
\nonumber &=&-\frac{1}{a\bigl(\ndivk\bigr)}\frac{1}{\log 2} \biggl\{\intab{\frac{1}{2}}{1} \Bigl(X_{n-[2ks],n}-U\bigl(\frac{n}{2ks}\bigr)\Bigr)\, \frac{ds}{s}-\intab{\frac{1}{2k}}{\frac{1}{k}}\Bigl(X_{n,n}-U\bigl(\frac{n}{ks} \bigr)\Bigr)\, \frac{ds}{s}\\
\nonumber  & & \mbox{\hspace{6.0cm}}-q\bigl(\frac{n}{k} \bigr)\Bigl(\intab{\frac{1}{2k}}{1}\frac{U\bigl(\frac{n}{ks}\bigr)-U\bigl(\frac{n}{2ks}\bigr)}{q\bigl(\frac{n}{k}\bigr)}\,\frac{ds}{s}-\log 2\Bigr)\biggr\}\\
\label{Integrals} &=& -\frac{1}{\log 2} \bigl\{J_1(k,n)- J_2(k,n)\bigr\}+\frac{q\bigl(\frac{n}{k}\bigr)}{a\bigl(\frac{n}{k}\bigr)}\frac{1}{\log 2} J_3(k,n).
\end{eqnarray}
By mimicking the steps of progression from  \eqref{2in1} to \eqref{2in1Exp}, we obtain for the first integral above that
\begin{eqnarray*}
\sqrt{2k}\, J_1(k,n)&:=&\sqrt{2k}\intab{\frac{1}{2}}{1} \frac{X_{n-[2ks],n}-U\bigl(\frac{n}{2ks}\bigr)}{a\bigl(\frac{n}{k} \bigr)}\, \frac{ds}{s}\\
	&=& \intab{\frac{1}{2}}{1}\sqrt{2k}\,R_2(s)\, \frac{ds}{s}+ \intab{\frac{1}{2}}{1}\Bigl(\frac{a\bigl(\frac{n}{2ks} \bigr)}{a\bigl(\frac{n}{k} \bigr)}-1 \Bigr)\sqrt{2k}\,R_2(s)\, \frac{ds}{s} \, .
\end{eqnarray*}
Hence, Proposition \ref{PropRTheta} while assuming that $\sqrt{k}\, a(n/k)/q(n/k)=O(1)$  (by appointment of Remark \ref{Rem2}) and application of the uniform bounds in \eqref{AuxUnifBounds} with $a_0(t):= a(t)\bigl(1+o(A(t)\bigr)$ and $A_0(t):= A(t)$, imply for each $\varepsilon >0$,
\begin{eqnarray*}
\sqrt{2k}\, J_1(k,n)&=&  \intab{\frac{1}{2}}{1} W_n(s)\bigl(1-\log(2 s)\bigr)\, \frac{ds}{s^2}+ o_p(1) \intab{\frac{1}{2}}{1}\log\bigl(\frac{1}{2s} \bigr)\bigl( \frac{1}{s}\bigr)^{3/2+\varepsilon}\, ds +o_p\Bigl(A_0\bigl(\frac{n}{k} \bigr) \Bigr)\,.
\end{eqnarray*}
Since the integral $\int_{1/2}^{1} W_n(s)(1-\log (2 s))s^{-2}\, ds$ converges to a sum of independent normal random variables, then the expression above allows to conclude that the first random component in \eqref{Integrals} is negligible with high probability because
\begin{equation*}
	J_1(k,n)=O_p\bigl(\frac{1}{\sqrt{k}}\bigr)\, .
\end{equation*}
Now, similarly to $I_2(k,n)$ in the proof of Proposition \ref{PropConsQ}, albeit under the second order condition \eqref{2ndPiofU} and pertaining uniform bounds provided by \eqref{2ndUnifBounds}, we now have that
\begin{eqnarray*}
J_2(k,n)&:=&\intab{\frac{1}{2k}}{\frac{1}{k}} \frac{X_{n,n}-U\bigl(\frac{n}{ks}\bigr)}{a\bigl(\frac{n}{k} \bigr)}\, \frac{ds}{s}\\
	&=& \frac{a(n)}{a\bigl(\ndivk\bigr)}\Bigl\{-\log 2\, \log \bigl(n U_{1,n}\bigr) + \frac{a_0(n)}{a(n)} \intab{\frac{1}{2}}{1} \log s\, \frac{ds}{s} + A_0(n) \intab{\frac{1}{2}}{1}\Bigl( \frac{(\log s)^2}{2}\pm \varepsilon s^{-\varepsilon}\Bigr)\, \frac{ds}{s} \Bigr\}.
\end{eqnarray*}
Again, note that $a_0(n)/a(n)-1= o\bigl(A(n)\bigr)$ and $A_0(n)= A(n)$. Hence,
\begin{eqnarray*}
\frac{a\bigl(\ndivk\bigr)}{a(n)}\frac{1}{\log 2}\,J_2(k,n)&=& -\log \bigl(n U_{1,n}\bigr)-\frac{\log 2}{2} +\frac{1}{\log 2}A(n) \intab{\frac{1}{2}}{1}\Bigl( \frac{(\log s)^2}{2}\pm \varepsilon s^{-\varepsilon}\Bigr)\, \frac{ds}{s} +o\bigl(A(n)\bigr)\\
	&=&  -\log \bigl(n U_{1,n}\bigr)-\frac{\log 2}{2} +o(1).
\end{eqnarray*}
Furthermore, assuming that $k=k(n)$ is such that $a(n)/a(n/k) \rightarrow 1$, then the following convergence in distribution holds
\begin{equation*}
\frac{a\bigl(\ndivk\bigr)}{a(n)}\frac{1}{\log 2}\,J_2(k,n) \conv{d} \Lambda -\frac{\log 2}{2},
\end{equation*}
where $\Lambda$ denotes a Gumbel random variable with distribution function $\exp\{-e^{-x}\},$ $x \in \real$ (cf. Eq. \eqref{Max} and subsequent text). The following also holds provided \eqref{AuxUnifBounds} and that $\sqrt{k}\, A(n/k)=O(1)$:
\begin{equation*}
\frac{1}{\log 2}\,J_2(k,n) \conv{d} \Lambda -\frac{\log 2}{2},
\end{equation*}

\noindent Finally we turn to the bias term $J_3(k,n)$. By assumption,
\begin{equation*}
 \frac{J_3(k,n)}{A(\tfrac{n}{k})} = \frac{1}{A(\tfrac{n}{k})}\biggl(\intab{\frac{1}{2k}}{1}\frac{U\bigl(\frac{n}{ks}\bigr)-U\bigl(\frac{n}{2ks}\bigr)}{q\bigl(\frac{n}{k}\bigr)}\,\frac{ds}{s}-\log 2\biggr) \,\, \arrowf{n} \,\,\lambda,
\end{equation*}
as $n\rightarrow \infty$. Therefore, since $A(\tfrac{n}{k}) \sim -a(\tfrac{n}{k})/q(\tfrac{n}{k})$ (cf. Remark \ref{Rem2}), the deterministic term $J_3(k,n)$ renders the following contribution to the asymptotic bias:
\begin{equation*}
\frac{q\bigl(\ndivk\bigr)}{a\bigl(\ndivk\bigr)}\frac{1}{\log 2}\,J_3(k,n) \arrowf{n} -\frac{\lambda}{\log 2}.
\end{equation*}
\end{pfofThm}

\begin{ex}\label{ExNegFrechRes}
We resume here the results for the Negative Fr\'echet distribution introduced in Example \ref{ExNegFrech}. The Negative Fr\'echet distribution with pertaining tail quantile function $U(t)= x^{F}-(\log t)^{-1/\beta}$, $t\geq 1$, $0<\beta<1$, satisfies the second order limiting condition \eqref{2ndExtRVIntU}  with $Q(t)= -(\beta\,\log t )^{-1}$.
\end{ex}

\bigskip

We are thus ready to pursue with devising the asymptotic distribution of $\hat{x}^F$. The following proposition rests heavily on the statement in Theorem \ref{TheoAN}.

\begin{prop}\label{PropANxF}
Under the conditions of Theorem \ref{TheoAN},
\begin{equation*}
\frac{1}{a(n/k)} \bigl(\hat{x}^F-x^F\bigr)-\frac{q\bigl(\ndivk\bigr)}{a\bigl(\ndivk\bigr)}\Bigl( \frac{\hat{q}\bigl(\ndivk\bigr)}{q\bigl(\ndivk\bigr)}-1\Bigr)\conv{P} 0.
\end{equation*}
\end{prop}
\begin{pf}
We use the fact that $X_{n-k,n}\id U(1/U_{k+1,n})$, where $U_{k+1,n}$ is the $(k+1)$th order statistic associated with a sample of $n$ independent and standard uniform random variables, in order to write
\begin{eqnarray*}
         \frac{\hat{x}^F-x^F}{a(n/k)}-\frac{q(n/k)}{a(n/k)}\Bigl(\frac{\hat{q}(n/k)}{q(n/k)}-1\Bigr)
                &=& \frac{\hat{x}^F-\hat{q}(n/k)}{a(n/k)}-\frac{x^{F}-q(n/k)}{a(n/k)}\\
                &=& \frac{X_{n-k,n}-U(n/k)}{a(n/k)}-\frac{U(\infty)-U(n/k)-q(n/k)}{a(n/k)}.
\end{eqnarray*}
Since $U\in \Pi(a)$ and $\sqrt{k}\bigl(k/(nU_{k+1,n})-1\bigr)$ is asymptotically standard normal (see Corollary 2.2.2 of \citet{deHaanFerreira:06}) then
\begin{equation*}
         \frac{X_{n-k,n}-U(n/k)}{a\bigl(\ndivk\bigr)}\, \id \, \frac{U\bigl(\frac{k}{nU_{k+1,n}}\ndivk\bigr)-U\bigl(\ndivk\bigr)}{a\bigl(\ndivk\bigr)}=O_p\Bigl(\frac{1}{\sqrt{k}}\Bigr)=o_p(1).
\end{equation*}
The rest follows from relation \eqref{MainRel}.
\end{pf}

The next theorem encloses an alternative formulation of the results comprised in Theorem \ref{TheoAN} and Proposition \ref{PropANxF} aiming at providing confidence bands for $\hat{x}^F$.

\begin{thm}\label{TheoANxF}
Let $X_1,\,X_2,\ldots$ be i.i.d. random variables with tail quantile function $U$ satisfying the second order condition \eqref{2ndPiofU}. Let $\hat{a}(n/k)$ be a consistent estimator for $a(n/k)$. Suppose  $k=k(n)$ is a sequence of positive integers such that, as $n\rightarrow \infty$, $k(n)\rightarrow \infty$, $k(n)/n\rightarrow 0$, $a(n)/a(n/k)\rightarrow 1$ and $\sqrt{k}\,A(n/k)=O(1)$. Furthermore assume that
\begin{equation*}
\limit{n}\, \frac{1}{A(n/k)}\biggl( \intab{\frac{1}{2k}}{1}\frac{U\bigl(\frac{n}{ks}\bigr)-U\bigl(\frac{n}{2ks}\bigr)}{q\bigl(\ndivk\bigr)}\, \frac{ds}{s}-\log 2\biggr)=\lambda\in \real.
\end{equation*}
Then
\begin{equation*}
\frac{1}{\hat{a}(n/k)}\bigl(\hat{x}^F-x^F\bigr) \conv{d} \Lambda -\frac{\log 2}{2}-\frac{\lambda}{\log 2}.
\end{equation*}
\end{thm}
\begin{pf}
The result follows easily by conjugating Theorem \eqref{TheoAN} with Proposition \eqref{PropANxF} and then applying Slustky's theorem.
\end{pf}

There are in the literature several possibilities for estimating the auxiliary (or scale) function $a(n/k)$. The most obvious choice is the Maximum Likelihood Estimator (MLE) by pretending that the exceedances over a certain high (random) threshold follow a Generalized Pareto distribution \citep[cf. section 3.4 of][]{deHaanFerreira:06}:
\begin{equation*}
	\hat{a}\bigl(\ndivk\bigr)=\hat{\sigma}^{MLE}:= \frac{1}{k}\sumab{i=0}{k-1} \bigl(X_{n-i,n}-X_{n-k,n}\bigr).
\end{equation*}

\section{Simulations}
\label{SecSims}

The three distributions intervening in this simulation study are taken throughout as key examples for the purpose of illustrating the finite sample behavior of our estimator for $x^F$ defined in \eqref{EstxF}.
\begin{description}
\item[Model 1:] Negative Fr\'echet, with distribution function $F(x)= 1-\exp\{-(x^F-x)^{-\beta}\}$, $x\leq x^F$, $\beta>0$. The pertaining tail quantile function $U$ is given by $U(t)= x^F-(\log t)^{-1/\beta}$, $t\geq 1$. Clearly $U\in \Pi(a)$ with $a(t)= \beta^{-1}(\log t)^{-1/\beta-1}$, $\beta >0$ (cf. Example \ref{ExNegFrech}).
\item[Model 2:] The distribution function $F$ given by $F(x)= 1-\exp\{-\tan(x/\beta)\}$, $0\leq x<\beta \pi /2$, $\beta >0$. The pertaining function $U$ is given by $U(t)= \beta \arctan (\log t)$, $t\geq 1$ and $U\in \Pi(a)$ with auxiliary function $a(t)= 1/\bigl( \log^2 t +\beta^{-2}\bigr)$.
\item[Model 3:] The distribution function $F$ given by $F(x)= 1-\exp\{(\pi/2)^{-\beta}-\bigl(\arcsin(1-x/\beta)\bigr)^{-\beta}\}$, $0\leq x< \beta$, $\beta>0$. The pertaining function $U$ is given by $U(t)= \beta\bigl\{1-\sin\bigl(\bigl[(2/\pi)^{\beta}+\log t\bigr]^{-1/\beta}\bigr)\}$, $t\geq 1$. Then $U\in \Pi(a)$ with $a(t)= (\log t)^{-(1/\beta+1)}\cos\bigl((\log t)^{-1/\beta}\bigr) $.

\end{description}

\begin{figure}
\caption{\footnotesize Probability density functions of Model 1 (\emph{first row}), Model 2 (\emph{second row}) and Model 3 (\emph{third row}).}
\begin{center}
\includegraphics[scale=0.3]{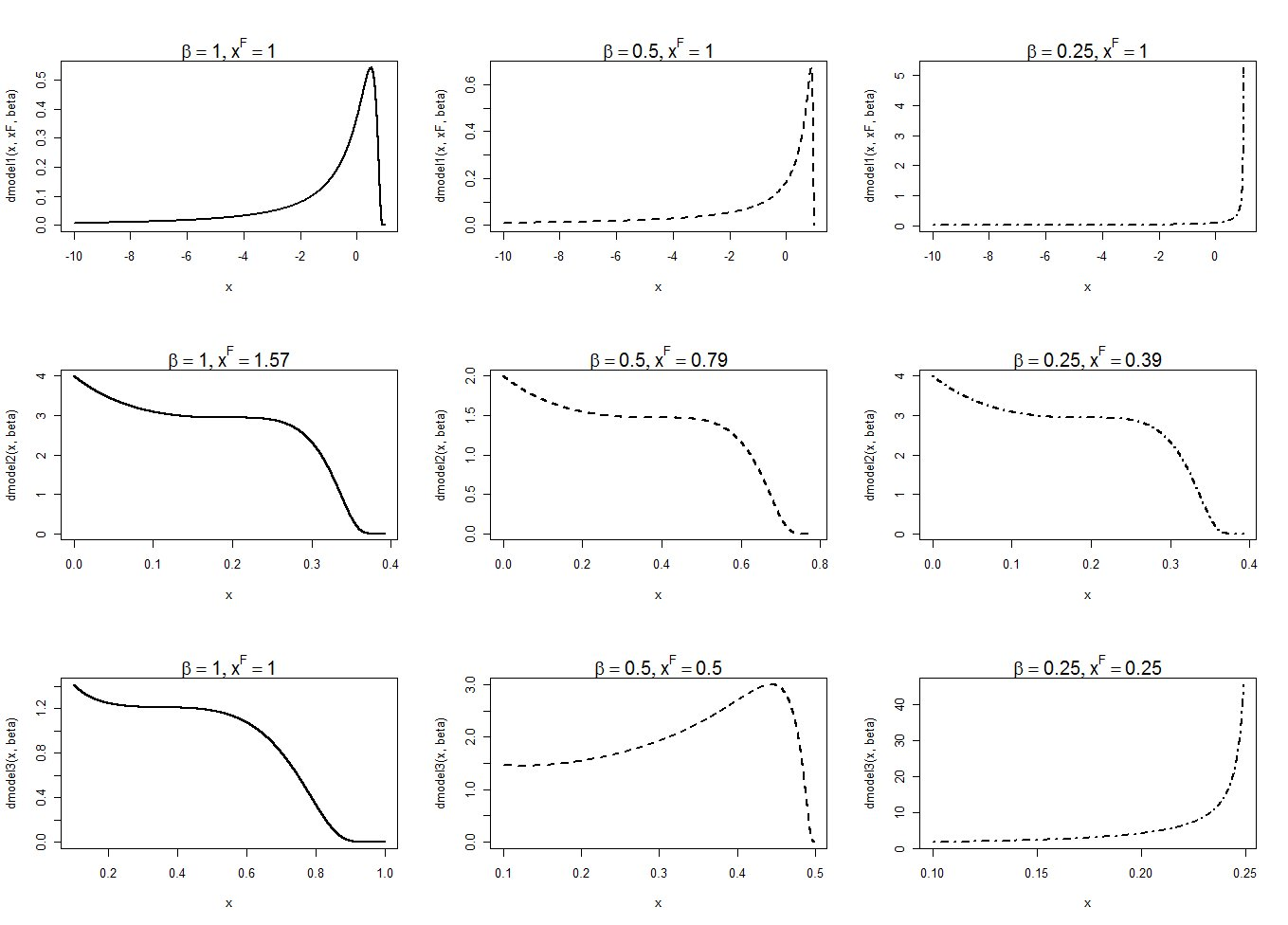}\hfill
\end{center}
\label{FigModels_pdfs}
\end{figure}

\begin{figure}
\caption{\footnotesize  Mean estimate and empirical Mean Squared Error  of  $\hat{x}^{F}$  for Model 1 with the true value $x^F=1$ and  several sample sizes: $n=100$ \emph{(first row)},  $n=1000$ \emph{(second row)},  $n=10000$ \emph{(third row)}; All plots are depicted against the number $k^*=2k$ of top order statistics used in the estimator.}
\begin{center}
\includegraphics[scale=0.22]{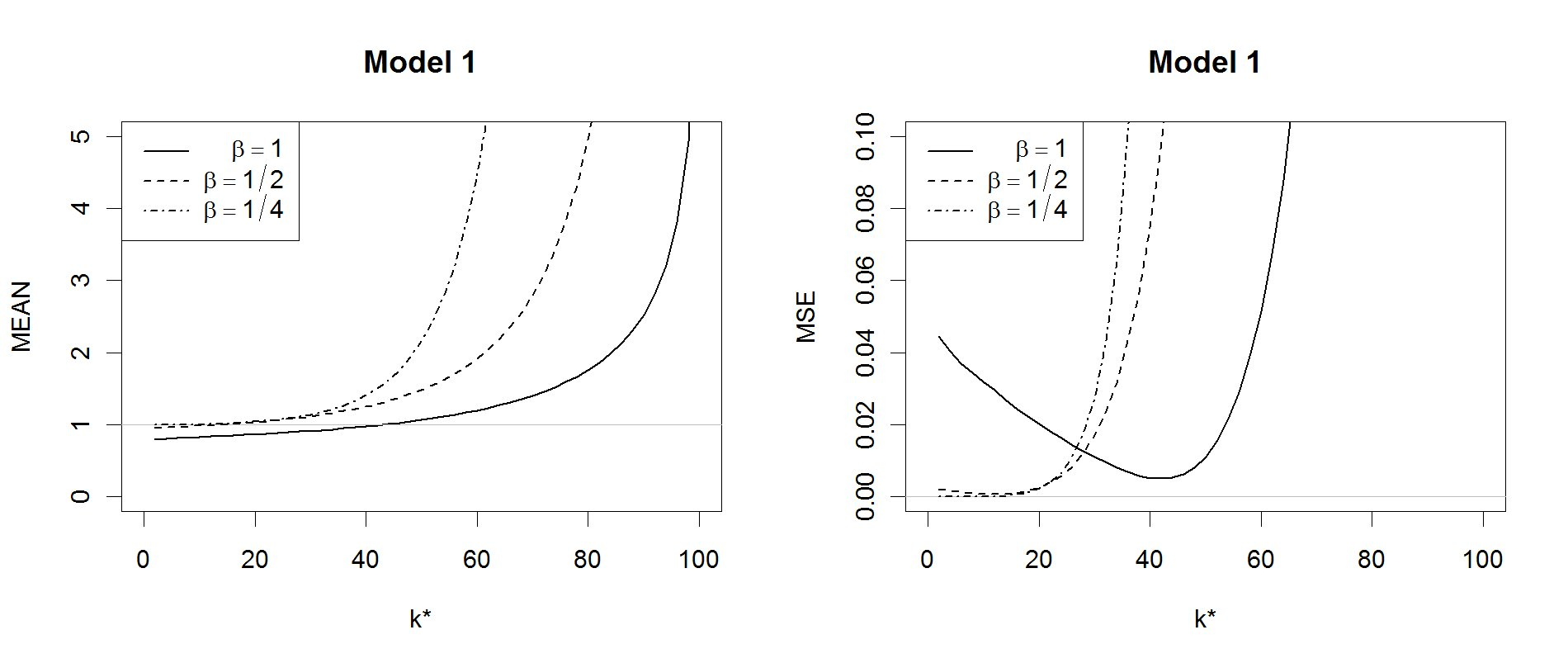}\hfill \includegraphics[scale=0.22]{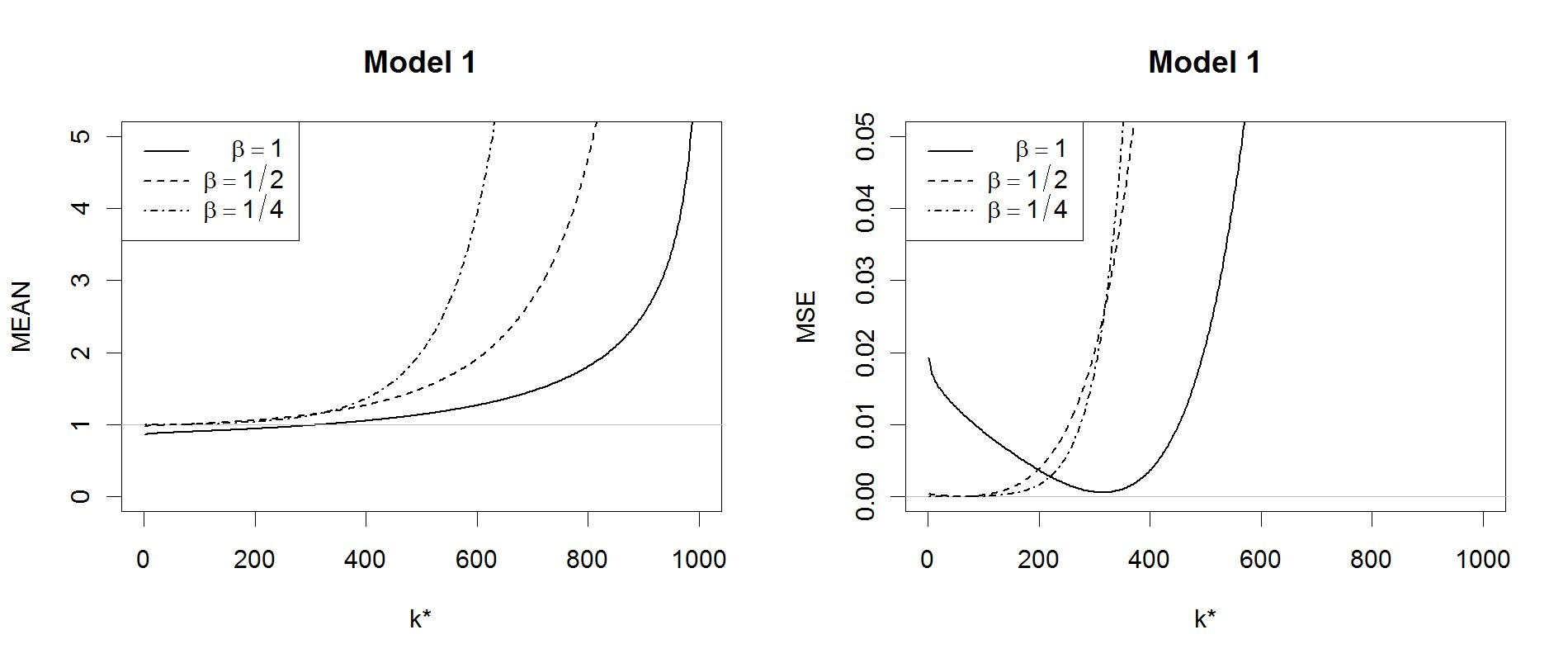}\hfill\includegraphics[scale=0.22]{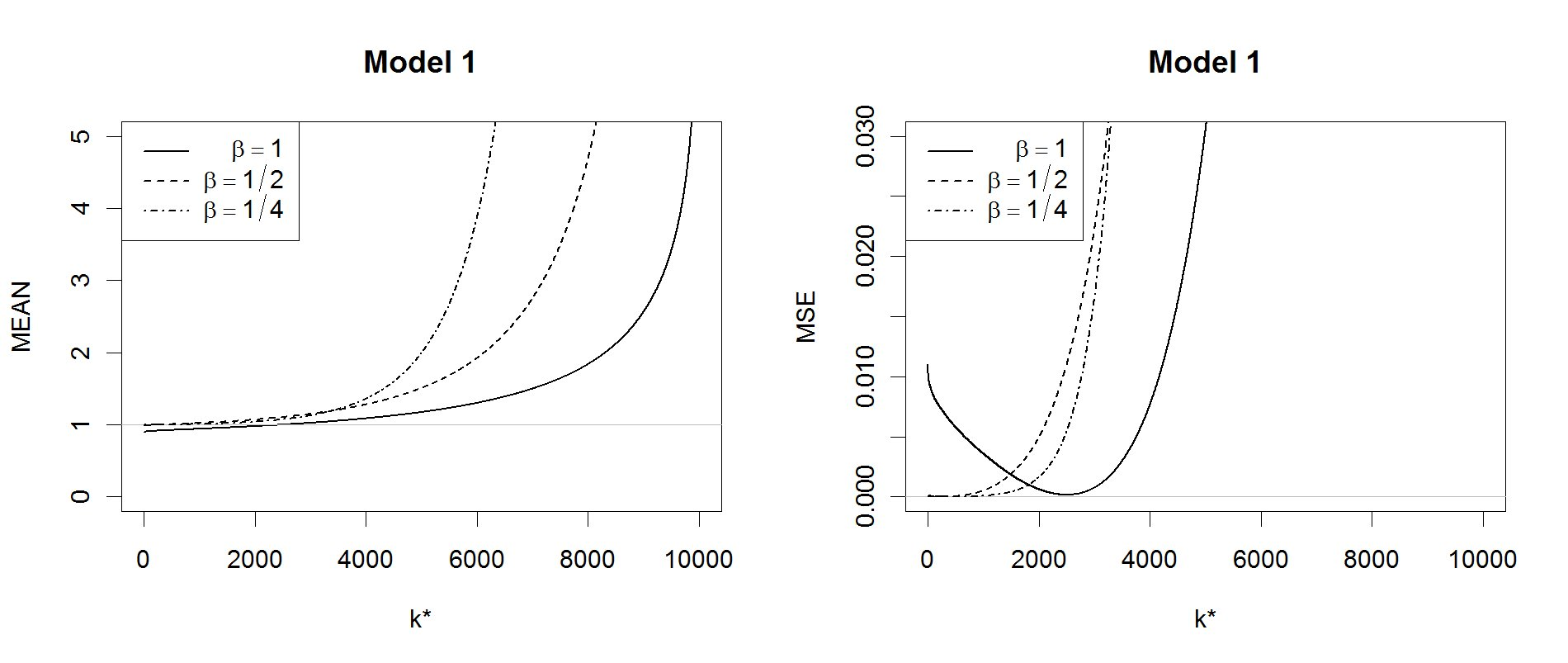}
\end{center}
\label{FigModel1}
\end{figure}

\begin{figure}
\caption{\footnotesize  Mean estimate and empirical Mean Squared Error  of  $\hat{x}^{F}$  for Model 2 with the true values  $x^F=\pi/8, \pi/4,\pi/2$  and  several sample sizes: $n=100$ \emph{(first row)},  $n=1000$ \emph{(second row)},  $n=10000$ \emph{(third row)}; All plots are depicted against the number $k^*=2k$ of top order statistics used in the estimator.}
\begin{center}
\includegraphics[scale=0.22]{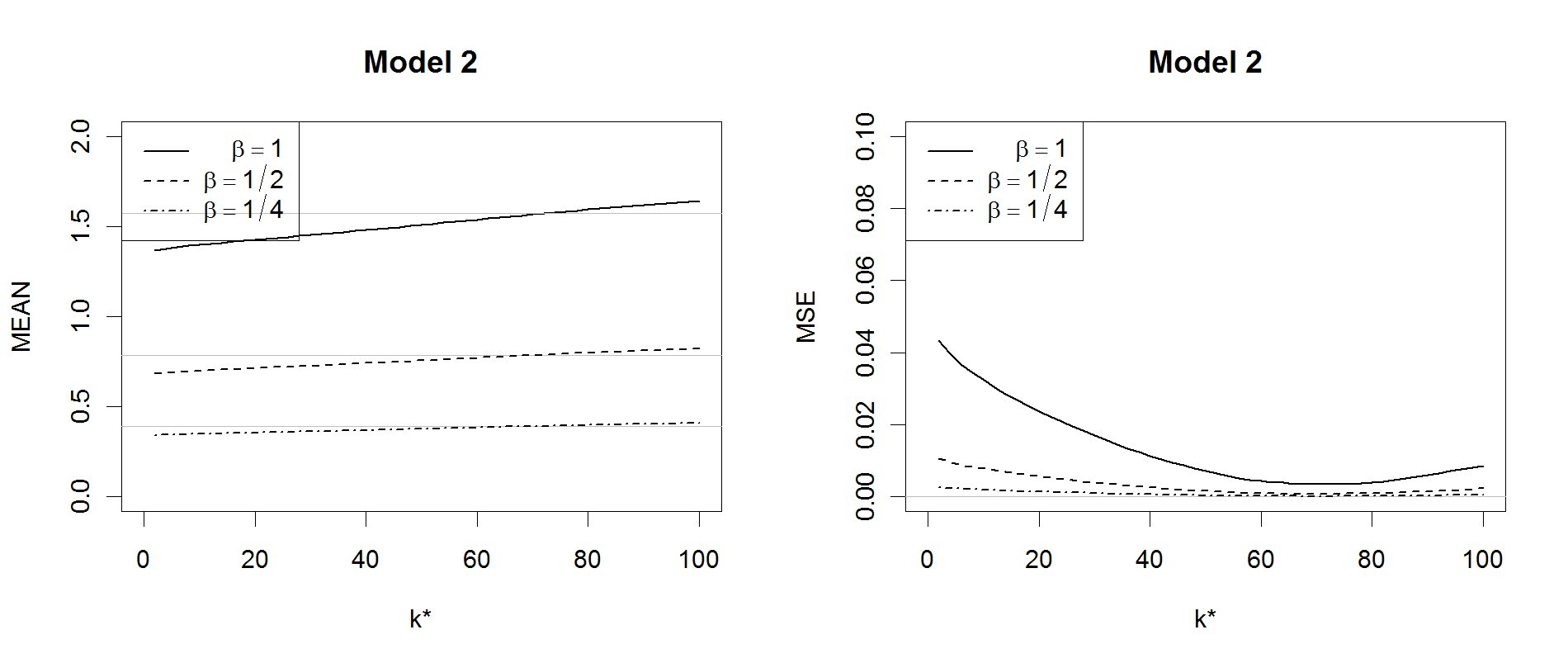}\hfill \includegraphics[scale=0.22]{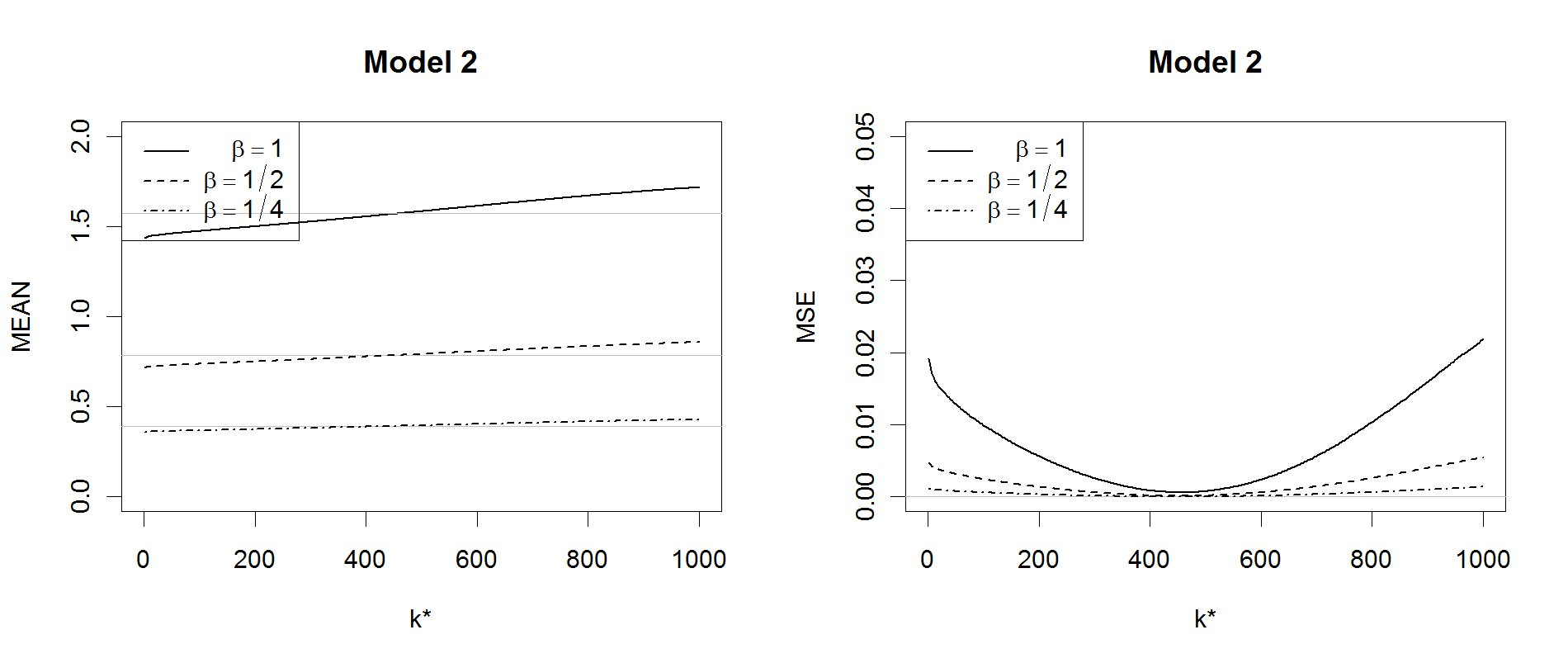}\hfill\includegraphics[scale=0.22]{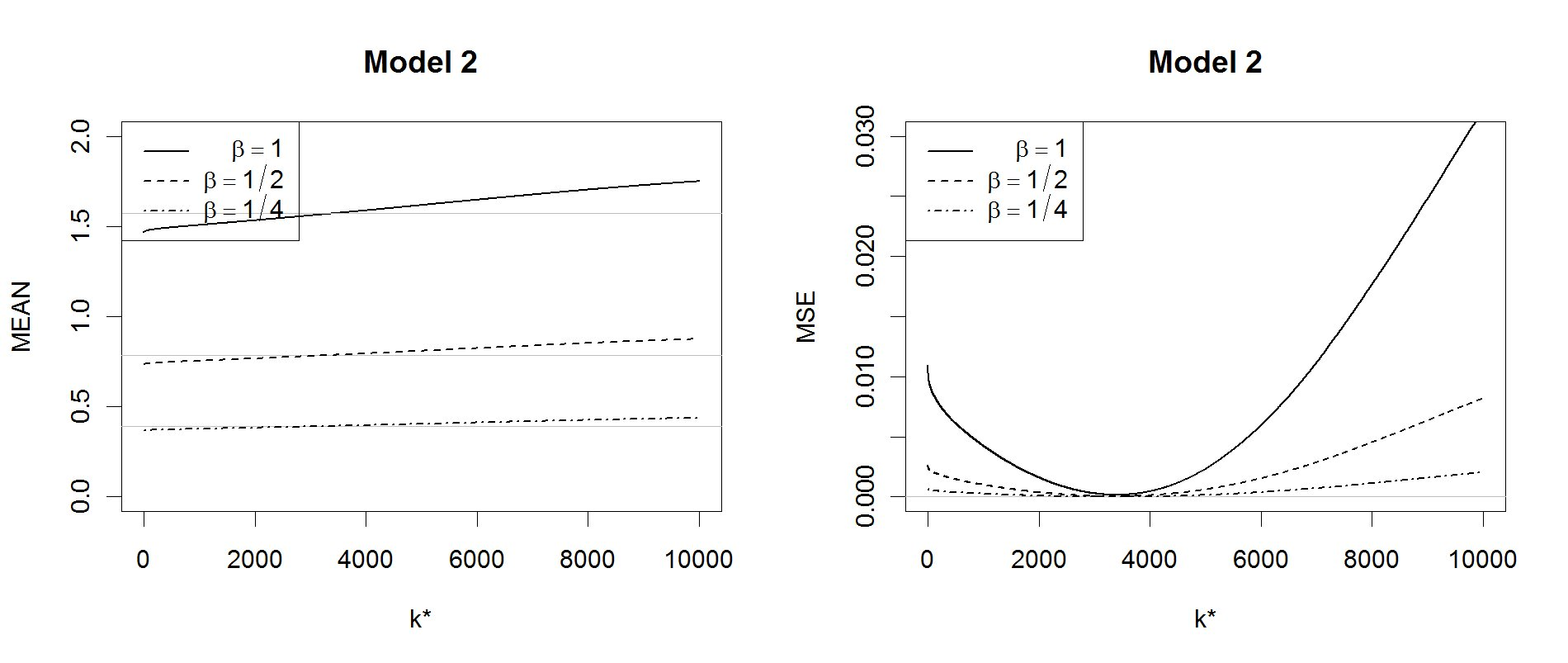}
\end{center}
\label{FigModel2}
\end{figure}

\begin{figure}
\caption{\footnotesize  Mean estimate and empirical Mean Squared Error  of  $\hat{x}^{F}$  for Model 3 with the true values $x^F=1/4, 1/2, 1$ and  several sample sizes: $n=100$ \emph{(first row)},  $n=1000$ \emph{(second row)},  $n=10000$ \emph{(third row)}; All plots are depicted against the number $k^*=2k$ of top order statistics used in the estimator.}
\begin{center}
\includegraphics[scale=0.22]{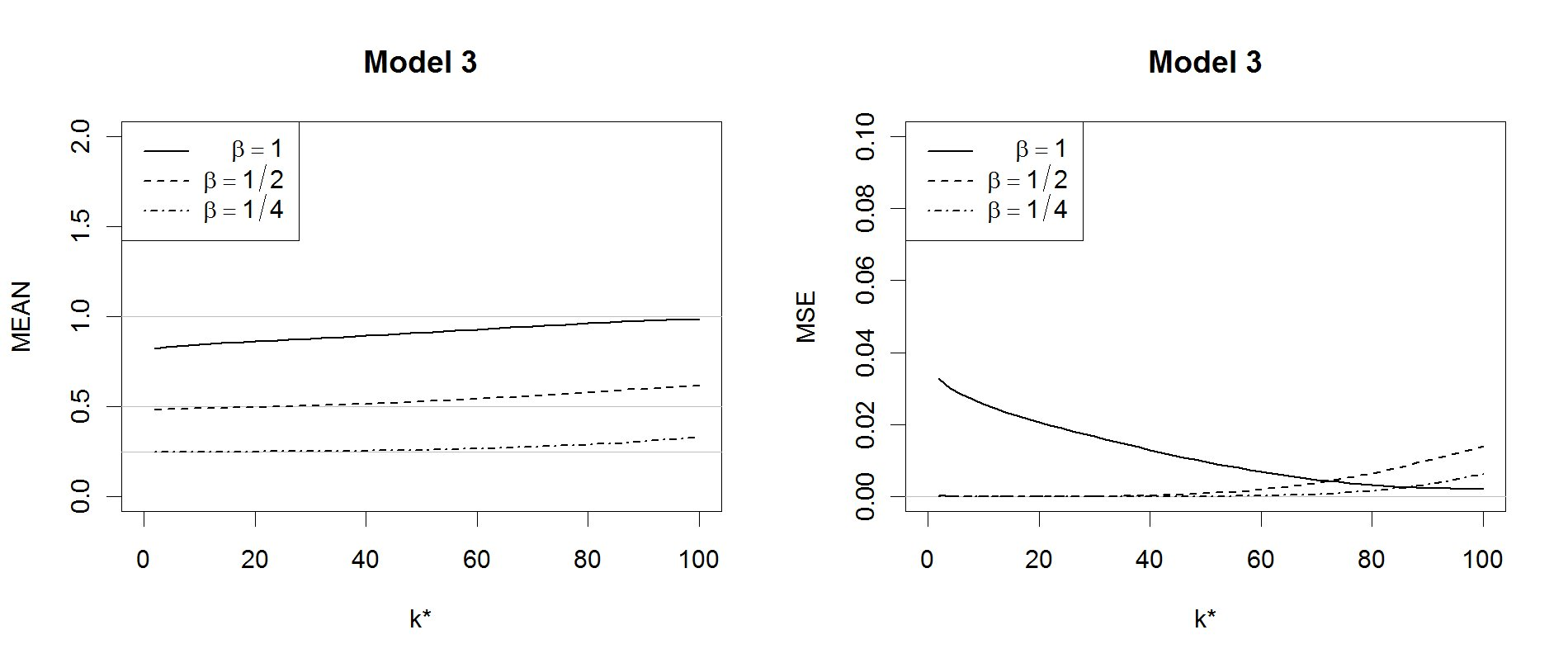}\hfill\includegraphics[scale=0.22]{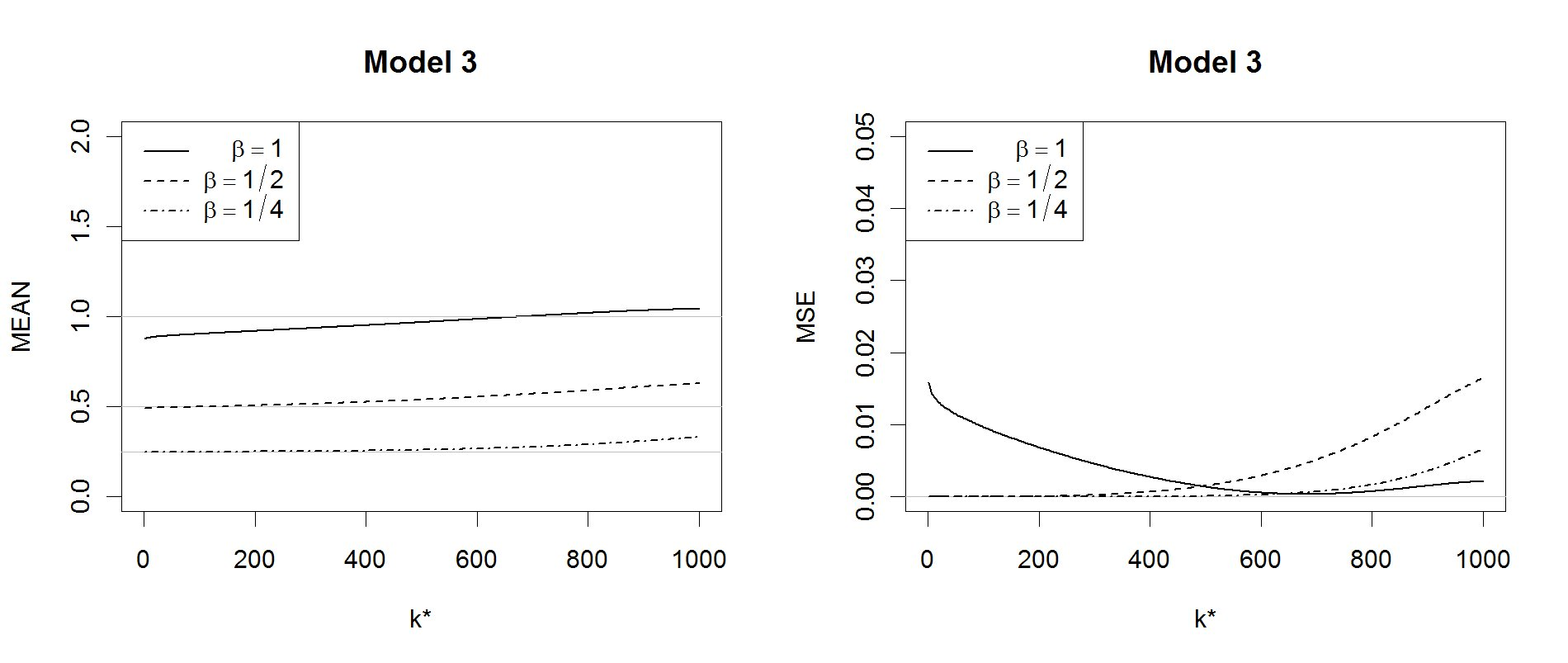}\hfill\includegraphics*[scale=0.22]{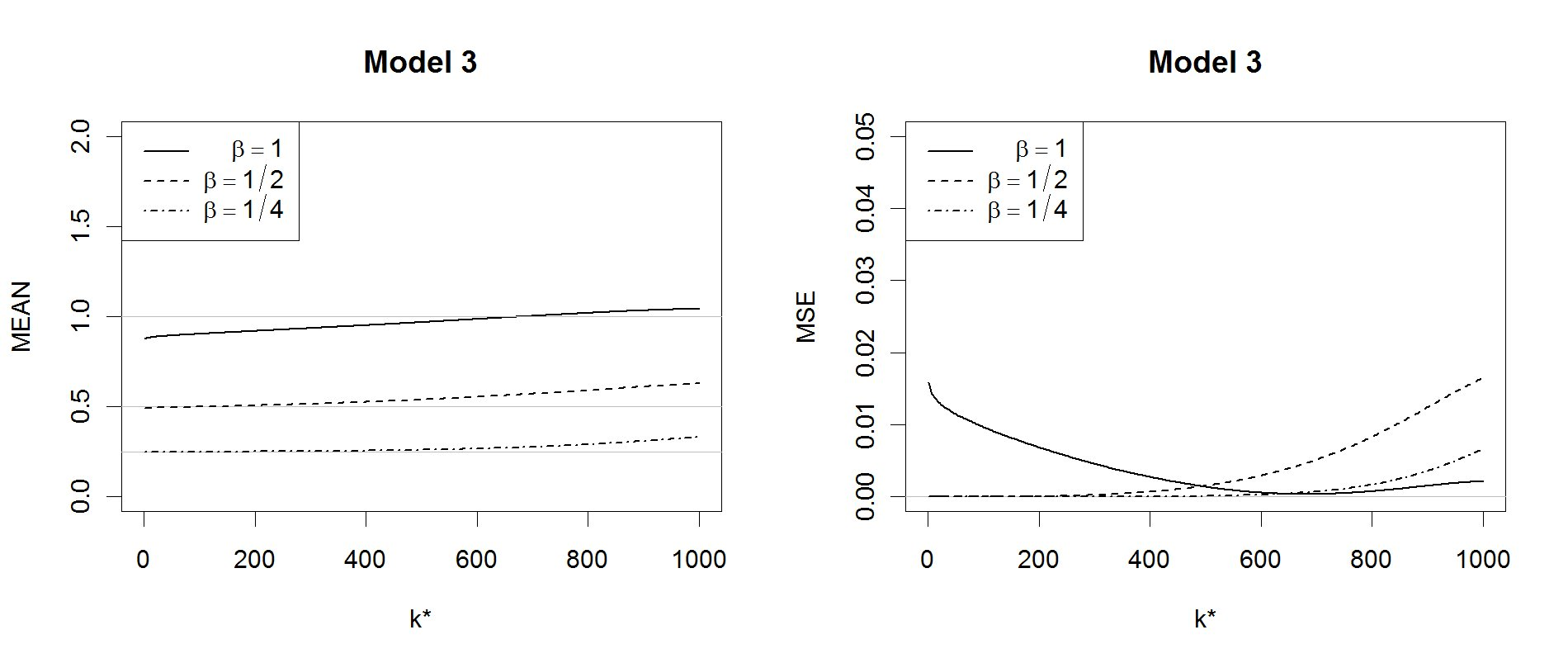}
\end{center}
\label{FigModel3}
\end{figure}

We have simulated $1000$ samples of size $n=100, 1000, 10000$, from each model and for different parameters $\beta=1, 1/2, 1/4$ . The results are depicted in Figures  \ref{FigModel1}, \ref{FigModel2} and \ref{FigModel3} . Since the number $k$ actually implies that the number of top order statistics used in the estimation is twice as much, we have plotted the estimated mean of $\hat{x}^F$ as a function of the latter, i.e., the plots are against $k^*=2k$. The most common approach of selecting the number $k$ (or $k^*$ in the present case) is to look for a region where the plots are relatively stable. This way, given the consistency property of the adopted estimator, one should in principle be away from small values of $k$ avoiding large variance (small $k$ is usually associated with a large variance) and not so far off in the tail preventing bias to instill (bias usually due to large $k$). As already discussed in Section \ref{SecAsympt}, for Model 1  an appropriate choice for an intermediate $k=k_n$ may be given by $k_n=(\log n)^r$, with $r\in(0,2]$. If we are using $n=1000$, for instance, and if we set  $r=2$,  the maximum allowed for $r$, we obtain $k\approx 48$ and thus $k^*\approx 96$. Bearing on  a value of $k^*$, around $100$ e.g., all the plots in Figure \ref{FigModel1} look quite stable in a close vicinity of the target value $x^F=1$ represented by the solid horizontal grey line.

A more thorough examination of the graphs in Figure  \ref{FigModel1}  seems to give accounts of a tendency to a better estimation under Model 1 (i.e., with underlying Negative Fr\'echet distribution) if the parameter $\beta$ is less than $1$, which corresponds to the case where the inherent second order conditions are satisfied. We recall that if $\beta \geq 1$, the Negative Fr\'echet distribution still satisfies the first order condition. Further details on the Negative Fr\'echet distribution are given in Examples \ref{ExNegFrech} and \ref{ExNegFrechRes}. Analogously,  in Figure  \ref{FigModel2} and   Figure  \ref{FigModel3}, the upper part of samples from Model 2 and Model 3 seems to yield small negative deviations from the true value $x^F$ specified in connection with the chosen values for the parameter $\beta>0$. However, the general pattern for these models is quite different in what concerns a  moderated  bias with increasing $k^*$, contrasting with the fast increasing bias with $k^*$ observed in model 1.
Note that for any model with right endpoint finite  the sample path of  $\hat{x}^F$ departures from the top value $x_{n,n}$, \ie, the sample maximum. 
    
    Taking all into account, we may conclude that the proposed estimator $\hat{x}^F$ performs reasonably well for parent distributions in the Gumbel domain detaining finite right endpoint $x^F$.

As a short final remark about the robustness of endpoint estimator defined in \eqref{Combi}, we can say it constitutes an advised inference procedure under Weibull domain of attraction. The theoretical background supporting this statement is a topic of further undergoing research, but  beyond the scope of the present subject.

\section*{Acknowledgement}

The authors are grateful to Professor Laurens de Haan for introducing
the appropriate characterization of distributions with finite right
endpoint in the Gumbel domain, at the origin of the proposed
estimator.


\bibliography{bibEndPoint}
\bibliographystyle{apalike}

\end{document}